\documentclass[10pt]{article}%

\usepackage{marvosym}

\usepackage{tikz}
\usepackage{pgf}
\usetikzlibrary{arrows}
\usetikzlibrary{calc}
\usetikzlibrary{decorations.pathmorphing}

\usepackage{soul}
\setstcolor{red}
\setulcolor{red}

\usepackage{rotating} 

\usepackage{mathrsfs}
\DeclareMathAlphabet{\mathpzc}{OT1}{pzc}{m}{it}

\usepackage{color}
\usepackage{amsmath}
\usepackage{graphicx}
\usepackage{amsfonts}
\usepackage{amssymb}%
\usepackage{latexsym}
\usepackage{psfrag}
\usepackage{accents}

\newtheorem{theorem}{Theorem}[section]
\newtheorem{corollary}[theorem]{Corollary}
\newtheorem{definition}[theorem]{Definition}
\newenvironment{proof}[1][Proof]{\noindent \emph{#1.} }
{\hfill \ \rule{0.5em}{0.5em}}
\newtheorem{lemma}[theorem]{Lemma}
\newtheorem{proposition}[theorem]{Proposition}

\newtheorem{assumption}[theorem]{Assumption}

\numberwithin{equation}{section}
\numberwithin{table}{section}
\numberwithin{figure}{section}

\usepackage[a4paper]{geometry}
\newtheorem{remark}[theorem]{Remark}
\newtheorem{example}[theorem]{Example}
\newcommand{\noi}{\noindent}

\newcommand{\cS}{{\cal S}}

\newcommand{\bx}{x}

\newcommand{\bn}{n}

\newcommand{\tildu}{\widetilde{u}}

\newcommand{\ri}{{\rm i}}
\newcommand{\rd}{{\rm d}}

\newcommand{\HoDk}{H^1_k(B_R)}

\newcommand{\beq}{\begin{equation}}
\newcommand{\eeq}{\end{equation}}
\newcommand{\beqs}{\begin{equation*}}
\newcommand{\eeqs}{\end{equation*}}
\newcommand{\bit}{\begin{itemize}}
\newcommand{\eit}{\end{itemize}}
\newcommand{\ben}{\begin{enumerate}}
\newcommand{\een}{\end{enumerate}}
\newcommand{\bal}{\begin{align}}
\newcommand{\eal}{\end{align}}
\newcommand{\bals}{\begin{align*}}
\newcommand{\eals}{\end{align*}}
\newcommand{\bse}{\begin{subequations}}
\newcommand{\ese}{\end{subequations}}
\newcommand{\bpr}{\begin{proposition}}
\newcommand{\epr}{\end{proposition}}
\newcommand{\bre}{\begin{remark}}
\newcommand{\ere}{\end{remark}}
\newcommand{\bpf}{\begin{proof}}
\newcommand{\epf}{\end{proof}}
\newcommand{\ble}{\begin{lemma}}
\newcommand{\ele}{\end{lemma}}
\newcommand{\bco}{\begin{corollary}}
\newcommand{\eco}{\end{corollary}}
\newcommand{\bex}{\begin{example}}
\newcommand{\eex}{\end{example}}
\newcommand{\bth}{\begin{theorem}}
\newcommand{\enth}{\end{theorem}}

\newcommand{\Rea}{\mathbb{R}}

\newcommand{\GR}{{\partial B_R}}

\newcommand{\eps}{\varepsilon}

\newcommand{\pdiff}[2]{\frac{\partial #1}{\partial #2}}

\newcommand{\gu}{\nabla u}

\newcommand{\vb}{\overline{v}}
\newcommand{\gvb}{\overline{\nabla v}}

\newcommand{\tendi}{\rightarrow \infty}

\def\XXint#1#2#3{{\setbox0=\hbox{$#1{#2#3}{\int}$}
     \vcenter{\hbox{$#2#3$}}\kern-.5\wd0}}

\usepackage{hyperref}
\definecolor{myblue}{rgb}{0,0,0.6}
\hypersetup{colorlinks=true,
linkcolor=myblue,citecolor=myblue,filecolor=myblue,urlcolor=myblue}
\newcommand*{\N}[1]{\left\|#1\right\|}

\allowdisplaybreaks[4]

\newcommand{\tfa}{\text{ for all }}
\newcommand{\tfor}{\text{ for }}

\newcommand{\tin}{\text{ in }}

\newcommand{\tand}{\text{ and }}
\newcommand{\tst}{\text{ such that }}
\newcommand{\tfind}{\text{ find }}

\newcommand{\vertiii}[1]{{\left\vert\kern-0.25ex\left\vert\kern-0.25ex\left\vert #1
    \right\vert\kern-0.25ex\right\vert\kern-0.25ex\right\vert}}

\newcommand{\DtN}{{\rm DtN}_k}

\usepackage{soul}

\definecolor{jwcol}{RGB}{27, 137, 18}  

\definecolor{dalcol}{rgb}{0.8,0,0}

\definecolor{escol}{rgb}{0,0,0.8}
\definecolor{estcol}{rgb}{0,0.5,0}
\definecolor{esnewcol}{rgb}{0,0.5,0}

\newcommand{\hFEM}{h}

\newcommand{\supp}{{\rm supp}}
\newcommand{\comp}{{\rm comp}}

\newcommand{\Cosc}{C_{\rm{osc}}}

\newcommand{\RR}{\mathbb{R}}

\newcommand{\OR}{B_R}

\newcommand{\truncbound}{\partial B_R}

\newcommand{\mymatrix}[1]{\mathsf{#1}}
\newcommand{\MA}{{\mymatrix{A}}}
\newcommand{\MI}{{\mymatrix{I}}}
\newcommand{\SPD}{{\mathsf{SPD}}}

\newcommand{\hatx}{\widehat{\bx}}

\newcommand{\Cqo}{{C_{\rm qo}}}

\newcommand{\Csol}{{C_{\rm sol}}}

\newcommand{\Ccont}{{C_{\rm cont}}}

\newcommand{\CDTN}{{C_{\rm DtN}}}

\newcommand{\mythmname}[1]{\textbf{\emph{(#1.)}}}

\newcommand{\uhigh}{u_{H^2}}
\newcommand{\ulow}{u_{\mathcal A}}
\newcommand{\vhigh}{v_{H^2}}
\newcommand{\vlow}{v_{\mathcal A}}

\newcommand{\Pilow}{\Pi_{L}}
\newcommand{\Pihigh}{\Pi_{H}}

\newcommand{\hsc}{\hbar}
\newcommand{\Psih}{\Psi_{\hsc}}
\newcommand{\WFh}{\operatorname{WF}_{\hsc}}
\newcommand{\pa}{\partial}

\newcommand{\Op}{{\rm Op}}

\DeclareMathOperator{\RC}{\mathsf{RC}}

\begin{document}

\title{Wavenumber-explicit convergence of the $hp$-FEM for the full-space heterogeneous Helmholtz equation with smooth coefficients}

\author{D.~Lafontaine\footnotemark[1]\,\,, E.~A.~Spence\footnotemark[2]\,\,, J.~Wunsch\footnotemark[3]}

\date{\today}

\footnotetext[1]{Department of Mathematical Sciences, University of Bath, Bath, BA2 7AY, UK, \tt D.Lafontaine@bath.ac.uk }
\footnotetext[2]{Department of Mathematical Sciences, University of Bath, Bath, BA2 7AY, UK, \tt E.A.Spence@bath.ac.uk }
\footnotetext[3]{Department of Mathematics, Northwestern University, 2033 Sheridan Road, Evanston IL 60208-2730, US, \tt jwunsch@math.northwestern.edu}

\maketitle

\begin{abstract}
A convergence theory for the $hp$-FEM applied to a variety of constant-coefficient Helmholtz problems was pioneered in the papers \cite{MeSa:10}, \cite{MeSa:11}, \cite{EsMe:12}, \cite{MePaSa:13}. This theory shows that, if the solution operator is bounded polynomially in the wavenumber $k$, then the Galerkin method is quasioptimal provided that $hk/p \leq C_1$ and $p\geq C_2 \log k$, where $C_1$ is sufficiently small, $C_2$ is sufficiently large, and both are independent of $k,h,$ and $p$. The significance of this result is that if $hk/p= C_1$ and $p=C_2\log k$, 
then quasioptimality is achieved with the total number of degrees of freedom proportional to $k^d$; i.e., the $hp$-FEM does not suffer from the pollution effect.

This paper proves the analogous quasioptimality result for the heterogeneous (i.e.~variable-coefficient) Helmholtz equation, posed in $\Rea^d$, $d=2,3$, with the Sommerfeld radiation condition at infinity, and $C^\infty$ coefficients. We also prove a bound on the relative error of the Galerkin solution 
in the particular case of the plane-wave scattering problem.
These are the first ever results on the wavenumber-explicit convergence of the $hp$-FEM for the Helmholtz equation with variable coefficients.
\end{abstract}


\section{Introduction}\label{subsec:intro}

\subsection{Context}\label{sec:1.1}

Over the last 10 years, a wavenumber-explicit convergence theory for the $hp$-FEM applied to the Helmholtz equation 
\beq\label{eq:Helmholtz}
\Delta u +k^2 u =-f 
\eeq
was established in the papers \cite{MeSa:10}, \cite{MeSa:11}, \cite{EsMe:12}, \cite{MePaSa:13}. This theory is based on decomposing solutions of the Helmholtz equation into two components:
\bit
\item[(i)] an analytic component, satisfying bounds with the same $k$-dependence as those satisfied by the full Helmholtz solution, and 
\item[(ii)] a component with finite regularity, satisfying bounds with improved $k$-dependence compared to those satisfied by the full Helmholtz solution.
\eit
Such a decomposition was obtained for
\bit
\item the Helmholtz equation \eqref{eq:Helmholtz} posed in $\Rea^d$, $d=2,3$, with compactly-supported $f$, and 
with the Sommerfeld radiation condition
\beq\label{eq:src}
\pdiff{u}{r}(\bx) - \ri k u(\bx) = o \left( \frac{1}{r^{(d-1)/2}}\right)
\eeq
as $r:= |\bx|\tendi$, uniformly in $\hatx:= \bx/r$ \cite[Lemma 3.5]{MeSa:10},
\item the Helmholtz exterior Dirichlet problem where the obstacle has analytic boundary \cite[Theorem 4.20]{MeSa:11},
\item the Helmholtz interior impedance problem where the domain is either smooth ($d=2,3$)  \cite[Theorem 4.10]{MeSa:11}, \cite[Theorem 4.5]{MePaSa:13}, or polygonal \cite[Theorem 4.10]{MeSa:11}, \cite[Theorem 3.2]{EsMe:12}.
\eit
This decomposition was then used to prove 
quasioptimality of the $hp$-FEM applied to the standard Helmholtz variational formulation in \cite{MeSa:10}, \cite{MeSa:11}, \cite{EsMe:12}, and applied to a discontinuous Galerkin formulation in \cite{MePaSa:13}. Indeed,
for the standard variational formulation (defined for the full-space problem in Definition \ref{def:vf} below) applied to the boundary value problems above, 
if the solution operator of the problem is bounded polynomially in $k$ (see Definition \ref{def:polybound} below), then there exist $C_1, C_2,$ and $\Cqo$ (independent of $k,h$, and $p$) such that if 
\beq\label{eq:threshold}
\frac{hk}{p}\leq C_1 \quad\tand\quad p\geq C_2 \log k
\eeq
then the Galerkin solution $u_N$ exists, is unique, and satisfies
\beqs
\N{u-u_N}_{H^1_k}\leq \Cqo \min_{v_N\in V_N} \N{u-v_N}_{H^1_k},
\eeqs
where $V_N$ is the $hp$ approximation space and the norm $\|\cdot\|_{H^1_k}$ is the standard weighted $H^1$ norm (defined by \eqref{eq:1knorm} below). 
Since the total number of degrees of freedom of the approximation space is proportional to $(p/h)^d$, 
the significance of this result is that it shows there is a choice of $h$ and $p$ such that the Galerkin solution is quasioptimal, 
with quasioptimality constant (i.e.~$\Cqo$) independent of $k$, and
with the total number of degrees of freedom proportional to $k^d$; thus, with these choices of $k$ and $p$, the $hp$-FEM does not suffer from the pollution effect \cite{BaSa:00}.

Over the last few years, there has been increasing interest in the numerical analysis of the heterogeneous Helmholtz equation, i.e.~the Helmholtz equation with variable coefficients
\beq\label{eq:Helmholtzvar}
\nabla\cdot(\MA \nabla u ) + k^2 n u = -f; 
\eeq
see, e.g., \cite{Ch:16}, \cite{BaChGo:17}, \cite{ChNi:20}, \cite{GaMo:19}, \cite{Pe:20}, \cite{GrSa:20}, \cite{GaSpWu:20}, \cite{LaSpWu:19}, \cite{GoGrSp:20}. 
However there do not yet exist in the literature analogous results to those in \cite{MeSa:10}, \cite{MeSa:11}, \cite{EsMe:12}, \cite{MePaSa:13} for the variable-coefficient Helmholtz equation. 

\subsection{Informal statement and discussion of the main results}\label{sec:informal}

\paragraph{The main results.}
This paper considers the variable-coefficient Helmholtz equation \eqref{eq:Helmholtzvar} with $C^\infty$ coefficients posed in $\Rea^d$, $d=2,3,$ with the Sommerfeld radiation condition at infinity. We obtain analogous results to those obtained in  \cite{MeSa:10} for this scenario with constant coefficients. That is, we prove quasioptimality of the $hp$-FEM under the conditions \eqref{eq:threshold} and provided that the solution operator is polynomially bounded in $k$; see Theorem \ref{thm:main2} below.

We obtain this result by decomposing the solution $u$ to \eqref{eq:Helmholtzvar} into two components:
\beqs
u|_{B_R}= \uhigh+ \ulow
\eeqs
where $\uhigh \in H^2(B_R)$ and $\ulow$ is analytic in $B_R$, where $B_R$ denotes the ball of radius $R$ centred at the origin (and $R$ is arbitrary); see Theorem \ref{thm:main1} below.
This is exactly analogous to the decomposition obtained in \cite{MeSa:10}, except that now $u$ satisfies the variable-coefficient equation \eqref{eq:Helmholtzvar} instead of \eqref{eq:Helmholtz}.


\paragraph{Overview of the ideas behind the decomposition and subsequent bounds.}
The idea in \cite{MeSa:10} was to decompose the data $f$ in \eqref{eq:Helmholtz}
 into ``low-'' and ``high-'' frequency components, with $\ulow$ the Helmholtz solution for the low-frequency component of $f$ and $\uhigh$ the Helmholtz solution for the high-frequency component of $f$. The frequency cut-offs were defining using the indicator function
\beq\label{eq:MScutoff}
1_{B_{\lambda k}}(\zeta) := 
\begin{cases}
1 & \tfor |\zeta|\leq \lambda\, k,\\
0 & \tfor |\zeta|\geq \lambda\, k,
\end{cases}
\eeq
with $\lambda$ a free parameter 
(see \cite[Equation 3.31]{MeSa:10} and the surrounding text). 
In \cite{MeSa:10} the frequency cut-off \eqref{eq:MScutoff} was then used with (a) the expression for $u$ as a convolution of the fundamental solution and the data $f$, and (b) the fact that the fundamental solution is known explicitly when $\MA=\MI$ and $n=1$, to obtain 
the appropriate bounds on $\ulow$ and $\uhigh$ using explicit calculation.

In this paper we use the same idea as in \cite{MeSa:10} of decomposing into low- and high-frequency components, but apply frequency cut-offs to the solution $u$ as opposed to the data $f$.
Then, given any cut-off function that is zero for $|\zeta|\geq Ck$, bounding the corresponding low-frequency component $\ulow$ is relatively straightforward using 
basic properties of the Fourier-transform (namely the expression for the Fourier transform of a derivative and Parseval's theorem). Indeed, in Fourier space each derivative corresponds to a power of the Fourier variable $\zeta$, and the frequency cut-off means that $|\zeta|\leq C k$ for $\ulow$; i.e.~every derivative of $\ulow$ brings down a power of $k$ compared to $\ulow$ (see \S\ref{sec:lowproof} below).
The main difficulty therefore is showing that the high-frequency component $\uhigh$ satisfies a bound with one power of $k$ improvement over the bound satisfied by $u$.

The main idea of the present paper is that the high-frequency cut-off can be chosen so that the (scaled) Helmholtz operator 
\beq\label{eq:operator}
P_k:= - \big(k^{-2}\nabla\cdot(\MA \nabla \cdot ) +  n\big)
\eeq
is \emph{semiclassically elliptic} on the support of the high-frequency cut-off. 
 Furthermore, choosing the cut-off function to be smooth (as opposed to discontinuous, as in \eqref{eq:MScutoff}) then allows us to use basic facts about the ``nice'' behaviour of elliptic semiclassical pseudodifferential operators (namely, they are invertible up to a small error) to prove the required bound on $\uhigh$.
(Recall that semiclassical pseudodifferential operators are just pseudodifferential operators with a large/small parameter; in this case the large parameter is $k$.) 

We now discuss further the frequency cut-offs and the bound on $\uhigh$ via ellipticity.

\paragraph{The frequency cut-offs.}

In contrast to \eqref{eq:MScutoff}, we choose $\chi_\mu \in C^\infty_{\comp}(\Rea^d)$ such that 
\beq\label{eq:cutoffintro}
\chi_\mu\big(k^{-2}|\zeta|^2\big) = 
\begin{cases}
1 & \tfor |\zeta|\leq \sqrt{\mu}\, k,\\
0 & \tfor |\zeta|\geq \sqrt{2\mu}\, k,
\end{cases}
\eeq
where the parameter $\mu$ is chosen later in the argument.
With the Fourier transform and its inverse defined by
\beq\label{eq:FT}
\mathcal F\varphi(\zeta) := \int_{\mathbb R^d} \exp\big( -\ri x \cdot \zeta \big)
\varphi(x) \, \rd x
\,\,\tand\,\,
\mathcal F^{-1}\psi(x) := (2\pi)^{-d} \int_{\mathbb R^d} \exp\big( \ri x \cdot \zeta \big)
 \psi(\zeta)\, \rd \zeta,
\eeq
we define the low-frequency cut-off $\Pilow$ by 
\beq\label{eq:Pilow_intro}
 \Pilow v (x)  := \mathcal{F}^{-1}\Big(\chi_\mu \big(k^{-2} |\zeta|^2 \big) \, \mathcal F v(\zeta)\Big),
\eeq
and the high-frequency cut-off $\Pihigh$ by 
\beq\label{eq:Pihi_intro}
 \Pihigh v (x)  := \mathcal{F}^{-1}\Big(\big(1-\chi_\mu \big(k^{-2} |\zeta|^2 \big)\big) \, \mathcal F v(\zeta)\Big),
\eeq
so that $\Pilow+ \Pihigh= I$.
We let $\varphi \in C^\infty_c$ be equal to one on $B_{R+1}$ and vanish outside $B_{R+2}$, and then 
\beq\label{eq:definitions_intro}
\ulow := \Pilow (\varphi u)\big|_{B_R}
\quad\tand\quad
\uhigh := \Pihigh(\varphi u)\big|_{B_R}.
\eeq

\paragraph{The bound on the high-frequency component $\uhigh$ via ellipticity.}

Recall that a PDE is elliptic if its principal symbol is non-zero. The concept of ellipticity for semiclassical differential operators (or, more generally, 
 semiclassical pseudodifferential operators) is analogous, except that it now involves the \emph{semiclassical principal symbol} (see \eqref{eq:ellip} below).
The semiclassical principal symbol of $P_k$ \eqref{eq:operator} is 
\beq\label{eq:symbol_intro}
\langle\MA\xi,\xi\rangle -n,
\eeq
where $\langle\cdot,\cdot\rangle$ denotes the $\ell^2$ inner product and $\xi= k^{-1}\zeta$ (see \eqref{eq:symbolPDE} below and the surrounding text). 

If the parameter $\mu$ in the cut-off function $\chi_\mu$
\eqref{eq:cutoffintro} is chosen to be a certain function of $\MA$ and $n$ (see
\eqref{eq:mu0} below), then the symbol \eqref{eq:symbol_intro} is
bounded away from zero when $k^{-2}|\zeta|^2 \geq \mu$, i.e.~in the
region of Fourier space where $\Pihigh$ is non-zero; one therefore describes $P_k$ as ``microlocally elliptic'', where the adjective ``microlocal" indicates that we have ellipticity on just a region of phase space (rather than on all of phase space in the more familiar global ellipticity).

These ellipticity properties are then used with the standard
microlocal elliptic estimate for pseudodifferential operators,
appearing in the semiclassical setting in, e.g., \cite[Appendix
E]{DyZw:19}, and stated in this setting as Theorem \ref{thm:elliptic}
below. The whole point is that a semiclassical pseudodifferential
operator that is elliptic in some region of phase space can be
inverted (up to some small error) in that region, and the norm of the
inverse is bounded uniformly in the large parameter (here $k$) as long
as one uses weighted norms (analogous to the familiar $H^1_k$ norm
\eqref{eq:1knorm}).

The result is that $\uhigh$ satisfies a bound with one power of $k$
improvement over the bound satisfied by $u$ (compare \eqref{eq:uhigh}
and \eqref{eq:Csol}).  To give a simple illustration of how
ellipticity can give this improved $k$-dependence, we contrast the
solutions of \beqs P_k u:= -(\Delta + k^2)u = f \quad\tand\quad
\widetilde{P}_k v:= -(\Delta-k^2)v=f , \eeqs with both equations posed
in $\Rea^d$ with compactly-supported $f$, and with $u$ satisfying the
Sommerfeld radiation condition \eqref{eq:src} and $v$ satisfying
boundedness at infinity.  The $L^2\rightarrow L^2$ bounds that are
sharp in terms of $k$-dependence are \beqs \N{u}_{L^2(B_R)} \lesssim
k^{-1} \N{f}_{L^2(\Rea^d)} \quad\tand\quad \N{v}_{L^2(\Rea^d)}
\lesssim k^{-2} \N{f}_{L^2(\Rea^d)}, \eeqs with the former given by
Part (i) of Theorem \ref{thm:polybound}, and the latter following from
the Lax-Milgram theorem.  The operator $P_k$ is 
  not semiclassically elliptic on
all of phase space (its semiclassical principal symbol is $|\xi|^2-1$), whereas $\widetilde{P}_k$ is semiclassically elliptic on all of phase space
(its semiclassical principal symbol is $|\xi|^2+1$); we
therefore see that ellipticity has resulted in the solution operator
having improved $k$-dependence. The proof of the bound on $\uhigh$ is
more technical, but the idea -- that the improvement in $k$-dependence
comes from ellipticity -- is the same.

\paragraph{The assumption that the solution operator is polynomially bounded in $k$.}

We need to assume that the solution operator is polynomially bounded in $k$ (in sense of Definition \ref{def:polybound} below), \emph{both} in 
proving the bound on $\uhigh$, \emph{and} in proving quasi-optimality of the $hp$-FEM.

The $k$-dependence of the Helmholtz solution operator depends on whether the problem is \emph{trapping} or \emph{nontrapping}.
For the heterogeneous Helmholtz equation \eqref{eq:Helmholtzvar} posed in $\Rea^d$ (i.e.~with no obstacle), trapping can be created by the coefficients $\MA$ and $n$; see, e.g., \cite{Ra:71}.
If the problem is nontrapping, then the Helmholtz solution operator (measured in the natural norms) is bounded in $k$. However, under the strongest form of trapping, the Helmholtz solution operator can grow exponentially in $k$ \cite{Ra:71}. Nevertheless, it has recently been proved that, if a set of frequencies of arbitrarily small measure is excluded, then the solution operator is polynomially bounded under any type of trapping \cite{LaSpWu:20}. Therefore, the result that the $hp$-FEM is quasi-optimal holds for a wide class of Helmholtz problems; see Corollary \ref{cor:qo} below.

\paragraph{Why do we need $C^\infty$ coefficients?} As highlighted above, our proof of the decomposition relies on standard results about semiclassical pseudodifferential operators (recapped in \S\ref{sec:recap}).
These results are usually stated for $C^\infty$ symbols, and thus to fit into this framework $\MA$ and $n$ must be $C^\infty$. However, 
examining the results we use, we see that we only 
need the symbol of the PDE to be in $C^{L}$ where $L$ depends only on the dimension $d$ and on the exponent $M$ appearing in the assumption that the solution operator is polynomially bounded (see Definitions \ref{def:Csol} and \ref{def:polybound} below).
Therefore, while we consider $\MA, n \in C^\infty$ to easily use results about semiclassical pseudodifferential operators from \cite{Zw:12}, \cite[Appendix E]{DyZw:19}, our results hold for $\MA \in C^{L}$ and $n  \in C^{L}$,
where $L=L(d,M)$.

\paragraph{Extending the decomposition result to the solution of other PDEs.}

Our proof of the decomposition result only relies on the principal symbol of the differential operator
being bounded below at infinity (in the sense of \eqref{eq:generalise1} below). Therefore, the decomposition result Theorem \ref{thm:main1} is valid for a much larger class of PDEs (and indeed pseudodifferential operators) than  \eqref{eq:Helmholtzvar}; see
 Remark \ref{rem:generalise} below for more details. 

In the follow-up paper \cite{LaSpWu:21}, we use the ideas of the present paper combined with much more sophisticated tools of semiclassical and microlocal analysis (namely the black-box scattering framework of Sj\"ostrand--Zworski \cite{SjZw:91}, the Helffer--Sj\"ostrand functional calculus \cite{HeSj:89}, and 
 associated results by Helffer, Robert, and Sj\"ostrand  \cite{HeRo:83}, \cite{Rob87}, \cite{Sj:97}) to prove analogous decompositions for a wide variety of scattering problems (albeit with slightly weaker estimates on $\ulow$). In particular, the main result of the present paper, Theorem \ref{thm:main1}, is rederived in this more general context as \cite[Theorem 1.16]{LaSpWu:21}.

We also note that, as announced in the abstract \cite{Me:20}, Bernkopf, Chaumont--Frelet, and Melenk are also studying the question of $k$-explicit convergence of the $hp$-FEM for the Helmholtz equation with variable coefficients.

\paragraph{Outline of the paper.}
\S\ref{sec:form} gives the definitions of the boundary-value problem and the finite-element method. \S\ref{sec:results} states the main results. \S\ref{sec:recap} recaps results about semiclassical pseudodifferential operators, with \cite{Zw:12} and \cite[Appendix E]{DyZw:19} as the main references.
 \S\ref{sec:proof_PDE} proves the result about the decomposition $u|_{B_R}= \uhigh+ \ulow$ (Theorem \ref{thm:main1}). \S\ref{sec:proof_FEM} proves the result about quasioptimality of the $hp$-FEM (Theorem \ref{thm:main2}).

\section{Formulation of the problem}
\label{sec:form}

\subsection{The boundary value problem}

\begin{assumption}[Assumptions on the coefficients]\label{ass:1}
$\MA \in C^{\infty}(\Rea^d, \SPD)$
(where $\SPD$ is the set of $d\times d$ real, symmetric, positive-definite matrices)
is  such that $\supp(\MI- \MA)$ is compact in $\Rea^d$ and
there exist $0<A_{\min}\leq A_{\max}<\infty$ such that, in the sense of quadratic forms,
\beq\label{eq:Alimits}
 A_{\min}\leq \MA(\bx)\leq A_{\max}\quad\text{ for all }\bx \in \Rea^d.
\eeq
$n\in C^\infty(\Rea^d,\Rea)$ is such that $\supp(1-n)$ is compact in $\Rea^d$ and there exist $0<n_{\min}\leq n_{\max}<\infty$ such that
\beq\label{eq:nlimits}
n_{\min} \leq n(\bx)\leq n_{\max}\quad \text{ for all } \bx \in \Rea^d.
\eeq
\end{assumption}

Let $R>0$ be such that $\supp(\MI- \MA) \cup \supp(1-n)\Subset B_R$,
where $B_R$ denotes the ball of radius $R$ about the origin and $\Subset$ denotes compact containment. 
Let $\gamma$ and $\partial_{\bn}$ denote the Dirichlet and Neumann traces, respectively, on $\partial B_R$, where the normal vector for the Neumann trace points out of $B_R$.

Define $\DtN: H^{1/2}(\truncbound) \rightarrow H^{-1/2}(\truncbound)$ to be the Dirichlet-to-Neumann map for the equation $\Delta u+k^2 u=0$ posed in the exterior of $B_R$ with the Sommerfeld radiation condition \eqref{eq:src}.
The definition of $\DtN$ in terms of Hankel functions and polar coordinates (when $d=2$)/spherical polar coordinates (when $d=3$) is given in, e.g., \cite[Equations 3.7 and 3.10]{MeSa:10}. 

\begin{definition}[Heterogeneous Helmholtz Problem on $\Rea^d$]\label{def:vf}
Given $\MA$ and $n$ satisfying Assumption \ref{ass:1}, $R>0$ such that $\supp(\MI- \MA) \cup \supp(1-n)\Subset B_R$, $k>0$, and $F \in (H^1(B_R))^*$, 
$u\in H^1(B_R)$ satisfies the \emph{Heterogeneous Helmholtz Problem on $\Rea^d$} if $u$ satisfies the variational problem 
\beq\label{eq:vf}
\text{ find } u \in H^1(B_R) \tst \quad a(u,v)=F(v) \quad \tfa v\in H^1(B_R),
\eeq
where
\begin{align}\label{eq:sesqui}
a(u,v)&:= \int_{B_R} 
\Big((\MA \nabla u)\cdot\overline{\nabla v}
 - k^2 n u\overline{v}\Big) - \big\langle \DtN (\gamma u),\gamma v\big\rangle_{\GR},
 \end{align}
where $\langle\cdot,\cdot\rangle_{\GR}$ denotes the duality pairing on $\GR$ that is linear in the first argument and antilinear in the second.
\end{definition}

\ble[Helmholtz boundary value problems included in Definition \ref{def:vf}]\label{lem:BVPs}

\

(i) If 
\beq\label{eq:L2data}
F(v) := \int_{B_R} f \, \overline{v}
\eeq
with $f\in L^2(B_R)$, then the solution $u$ to \eqref{eq:vf} equals $\tildu |_{B_R}$, where $\tildu \in H^1_{\rm loc}(\Rea^d)$ is the solution to 
\beqs
\nabla\cdot(\MA \nabla \tildu ) + k^2 n \tildu = -f \quad\tin \Rea^d,
\eeqs
and $\tildu$ satisfies the Sommerfeld radiation condition \eqref{eq:src}. 

(ii) If 
\beq\label{eq:planewavedata}
F(v):= \int_{\truncbound}\left(\partial_{\bn}u^I - \DtN(\gamma u^I)\right)\overline{\gamma v}
\quad\text{ with } \quad u^I(\bx) := \exp(\ri k \bx\cdot a),
\eeq
where $a\in \Rea^d$ with $|a|=1$, 
then the solution $u$ to \eqref{eq:vf} equals
$\tildu |_{B_R}$, where $\tildu \in H^1_{\rm loc}(\Rea^d)$ is the solution of the \emph{Helmholtz plane-wave scattering problem}; i.e. 
\beqs
 \nabla\cdot(\MA\nabla \tildu)+k^2 n \tildu =0 \quad\tin\Rea^d,
\eeqs
and $\tildu^S:= \tildu-u^I$ satisfies the Sommerfeld radiation condition \eqref{eq:src}. 
\ele

\noi Part (i) of Lemma \ref{lem:BVPs} is proved in, e.g., \cite[Lemma 3.3]{GrPeSp:19}; the proof of Part (ii) is similar.

Let the weighted $H^1$ norm, $\|\cdot\|_{H^1_k(B_R)}$, be defined by
\begin{equation} \label{eq:1knorm}
\|u\|^2_{H^1_k(B_R)} :=\N{\nabla u}_{L^2(B_R)}^2 + k^2 \N{u}_{L^2(B_R)}^2.
\eeq

\ble
\label{lem:wellposed}
The solution of the Heterogeneous Helmholtz Problem on $\Rea^d$ (defined in Definition \ref{def:vf}) exists, is unique, and there exists $C(k,\MA,n, R)>0$
such that 
\beq\label{eq:Fredholm}
\N{u}_{\HoDk} \leq C \N{F}_{(\HoDk)^*} \quad\tfa k>0.
\eeq
\ele

\bpf
Uniqueness follows from the unique continuation principle; see \cite[\S1]{GrPeSp:19}, \cite[\S2]{GrSa:20} and the references therein.
Since $a(\cdot,\cdot)$ satisfies a G\aa rding inequality (see \eqref{eq:Garding} below), 
Fredholm theory then gives existence and the bound \eqref{eq:Fredholm}.
\epf

\paragraph{Properties of $\DtN$ and $a(\cdot,\cdot)$.}

We use later the following two properties of $\DtN$:~given $k_0, R_0>0$, there exists $\CDTN=\CDTN(k_0 R_0)$ such that, for all $k\geq k_0$ and $R\geq R_0$,
\beq\label{eq:CDTN1}
\big|\big\langle \DtN(\gamma u), \gamma v\rangle_{\partial B_R}\big\rangle\big| \leq \CDTN_1 \N{u}_{H^1_k(\OR)}  \N{v}_{H^1_k(\OR)} 
\eeq
for all $u,v \in H^1(\OR)$, 
and 
\beq\label{eq:CDTN2}
- \Re \big\langle \DtN \phi,\phi\big\rangle_{\GR} \geq 0
\quad\tfa \phi \in H^{1/2}(\GR).
\eeq
For a proof of \eqref{eq:CDTN1}, see \cite[Lemma 3.3]{MeSa:10}. For a proof of \eqref{eq:CDTN2}, see \cite[Theorem 2.6.4]{Ne:01} (for $d=3$) and \cite[Corollary 3.1]{ChMo:08} or \cite[Lemma 3.10]{MeSa:10} (for $d=2,3$).

Let $\Ccont = \Ccont(\MA, n, R, k_0)$ be the \emph{continuity constant} of the sesquilinear form $a(\cdot,\cdot)$ (defined in \eqref{eq:sesqui}) in the norm $\|\cdot\|_{H^1_k(B_R)}$; i.e.
\beqs
\big| a(u,v)\big|\leq \Ccont \N{u}_{H^1_k(B_R)} \N{v}_{H^1_k(B_R)} \quad\tfa u, v \in H^1(B_R) \text{ and } k\geq k_0.
\eeqs
By the Cauchy-Schwarz inequality and \eqref{eq:CDTN1},
\beq\label{eq:Cc}
\Ccont \leq \max\{ A_{\max}, n_{\max} \} + \CDTN_1.
\eeq

\subsection{The behaviour of the solution operator for large $k$}

\begin{definition}[$\Csol$]\label{def:Csol}
Given $f\in L^2(B_R)$, let $u$ be the solution of the heterogeneous Helmholtz equation \eqref{eq:Helmholtzvar} with the Sommerfeld radiation condition \eqref{eq:src} (i.e.~$u$ is the solution of the variational problem \eqref{eq:vf} with $F(v)$ given by \eqref{eq:L2data}).
Given $k_0>0$, let $\Csol= \Csol(k,\MA,n,R,k_0)>0$ be such that
\beq\label{eq:Csol}
\N{u}_{\HoDk} \leq \Csol \N{f}_{L^2(B_R)} \quad\tfa k>0.
\eeq
\end{definition}
$\Csol$ exists by Lemma \ref{lem:wellposed}; indeed, with $C$ given by \eqref{eq:Fredholm}, $\Csol:= C/k$.

How $\Csol$ depends on $k$ is crucial to the analysis below, and to emphasise this we write $\Csol=\Csol(k)$. Below we consider $\Csol$ with different values of $R$, and we then write, e.g., $\Csol(k;R)$ (as in the bound \eqref{eq:ulow} below).

A key assumption in the analysis of the Helmholtz $hp$-FEM is that $\Csol(k)$ is polynomially bounded in $k$ in the following sense.

\begin{definition}[$\Csol$ is polynomially bounded in $k$]\label{def:polybound}
Given $k_0$ and $K \subset [k_0,\infty)$, $\Csol(k)$ is polynomially bounded for $k\in K$ if 
 there exists $C>0$ and $M>0$ such that
\beq\label{eq:polybound}
\Csol(k) \leq C k^M \tfa k \in K,
\eeq
where $C$ and $M$ are independent of $k$ (but depend on $k_0$ and possibly also on $K, \MA, n, d, R$).
\end{definition}

There exist $C^\infty$ coefficients $\MA$ and $n$ such that $\Csol(k_j) \geq c_1 \exp (c_2 k_j)$ for $0<k_1<k_2<\ldots$ with $k_j\tendi$ as $j\tendi$, see \cite{Ra:71}, but this exponential growth is the worst-possible, since 
$\Csol(k) \leq c_3 \exp( c_4 k)$ for all $k\geq k_0$ by \cite[Theorem 2]{Bu:98}.
We now recall results on when $\Csol(k)$ is polynomially bounded in $k$.

\begin{theorem}[Conditions under which $\Csol(k)$ is polynomially bounded in $k$]\label{thm:polybound}

\

(i) $\MA$ and $n$ are $C^\infty$ and nontrapping (i.e.~all the trajectories of the Hamiltonian flow defined by the symbol of \eqref{eq:Helmholtzvar} starting in $B_R$ leave $B_R$ after a uniform time), then $\Csol(k)$ is independent of $k$ for all $k$, i.e., \eqref{eq:polybound} holds for all $k$ with $M=0$.

(ii) If $n=1$ and $\MA$ is $C^{0,1}$ then, given $k_0>0$ and $\delta>0$ there exists a set $J\subset [k_0,\infty)$ with $|J|\leq \delta$ such that  
\beq
\Csol(k) \leq C k^{5d/2+1+\eps} \quad\tfa k \in [k_0,\infty)\setminus J,
\eeq
 for any $\eps>0$, where $C$ depends on $\delta, \eps, d, k_0,$ and $\MA$.
If $\MA$ is $C^{1,\sigma}$ for some $\sigma>0$ then the exponent is reduced to $5d/2+\eps$.
\end{theorem}

\bpf[References for the proof]

(i) is proved using \emph{either} 
 (a) the propagation of singularities results of \cite{DuHo:72} combined with either the parametrix argument of \cite[Theorem 3]{Va:75}/ \cite[Chapter 10, Theorem 2]{Va:89} or Lax--Phillips theory \cite{LaxPhi}, \emph{or} (b) the defect-measure argument of \cite[Theorem 1.3 and \S3]{Bu:02}.
It has recently been proved that, for this situation, $\Csol$ is proportional to the length of the longest trajectory in $B_R$; see \cite[Theorems 1 and 2, and Equation 6.32]{GaSpWu:20}. 

(ii) is proved in \cite[Theorem 1.1 and Corollary 3.6]{LaSpWu:20}. 
\epf

\subsection{The finite-element method}\label{sec:FEMsetup}

Let $(V_N)_{N=0}^\infty$ be a sequence of finite-dimensional subspaces of $H^1(B_R)$ that converge to $H^1(B_R)$ in the sense that, for all $v\in H^1(B_R)$,
\beqs
\lim_{N\tendi} \Big(\min_{v_N \in V_N} \N{v-v_N}_{H^1(B_R)}\Big) =0.
\eeqs
Later we specialise to the triangulations described in \cite[\S5]{MeSa:10}, which allow curved elements and thus fit $\partial B_R$ exactly.

The finite-element method for the variational problem \eqref{eq:vf} is the Galerkin method applied to the variational problem \eqref{eq:vf}, i.e.
\beq\label{eq:FEM}
\text{ find } u_N \in V_N \tst\,\, a(u_N,v_N)=F(v_N) \,\, \tfa v_N\in  V_N.
\eeq

\section{Statement of the main results}\label{sec:results}

\begin{theorem}[Decomposition of the solution]\label{thm:main1}
Let $\MA$ and $n$ satisfy Assumption \ref{ass:1} and let $R>0$ be such that $\supp(\MI-\MA)\cup\supp(1-n)\Subset B_R$. Given $f\in L^2(B_R)$, let $u$ satisfy $\nabla\cdot(\MA \nabla u ) + k^2 n u = -f$ in $\Rea^d$ and the Sommerfeld radiation condition \eqref{eq:src}.

If $\Csol(k)$ is polynomially bounded (in the sense of Definition \ref{def:polybound}) for $k\in K\subset [k_0,\infty)$, then
there exist $C_3, C_4, C_5>0$ such that
\beqs
u|_{B_R}= \uhigh+ \ulow
\eeqs
where $\uhigh \in H^2(B_R)$ with
\beq\label{eq:uhigh}
\N{\partial^\alpha \uhigh}_{L^2(B_R)} \leq C_3 k^{|\alpha|-2} \N{f}_{L^2(B_R)}
 \quad\text{ for all $|\alpha|\leq 2$ and for all $k\in K\subset [k_0,\infty)$}, 
\eeq
and $\ulow \in C^\infty(B_R)$ with
\beq\label{eq:ulow}
\N{\partial^\beta \ulow}_{L^2(B_R)} \leq \Csol(k;R+2)\, \, C_4\, \big(C_5 k\big)^{|\beta|-1} \N{f}_{L^2(B_R)} \quad\text{ for all $\beta$ and for all $
k\in K\subset [k_0,\infty)$},
\eeq
where $C_3, C_4,$ and $C_5$ depend on $\MA, n, d$, and $k_0$, but are independent of $k$, $f$, $\alpha$, and $\beta$.
\end{theorem}

\bre[$\ulow$ is analytic]
Since $C_4$ and $C_5$ are independent of $\beta$, the bound \eqref{eq:ulow} implies that $\ulow$ is in the class of analytic functions on $B_R$,  $\mathcal{A}(B_R)$, defined by
\begin{align*}
\mathcal{A}(B_R) := \bigg\{ v \in \bigcap_{n\in \mathbb{N}} H^n(B_R) : \exists \,c_0, c_1 >0, \text{ independent of $n$,} \tst \,| u |_{H^n(B_R)}\leq c_1 c_0^n n! \bigg\},
\end{align*}
where $|u|_{H^n}^2:= \sum_{|\alpha|=n}\N{\partial^\alpha u}_{L^2}^2$.
See, e.g., \cite[\S1.1.b]{CoDaNi:10}, both for this definition, and for how the definition implies convergence of the Taylor series of elements of $\mathcal{A}(B_R)$ at every point in $\overline{B_R}$.
\ere

\bre[The bounds of Theorem \ref{thm:main1} written with the notation $\nabla^n$]\label{rem:nabla}
The analogous bounds to \eqref{eq:uhigh} and \eqref{eq:ulow} in \cite{MeSa:10}, \cite{MeSa:11} are written using the notation 
\beqs
\big| \nabla^n u(x)\big|^2 := \sum_{|\alpha|=n} \frac{n!}{\alpha!} \big| \partial^{\alpha} u(x)\big|^2.
\eeqs
Since $\sum_{|\alpha|=n} (n!/\alpha!) = d^n$,
\beqs
\text{ if }
\,\,
\N{\partial^\alpha u}_{L^2(B_R)} \leq \mathcal{C}_1 \, \big(\mathcal{C}_2\big)^{|\alpha|} \tfa \alpha \text{ with $|\alpha|=n$},\,\,
\text{ then } 
\N{\nabla^n u}_{L^2(B_R)} \leq \mathcal{C}_1 \, \big(\mathcal{C}_2 \sqrt{d} \big)^{n},
\eeqs
and so the bounds \eqref{eq:uhigh} and \eqref{eq:ulow} can also be written as bounds on $\N{\nabla^n \uhigh}_{L^2(B_R)}$ and $\N{\nabla^n \ulow}_{L^2(B_R)}$ respectively.
\ere

The following result about quasioptimality of the $hp$-FEM is then obtained by combining Theorem \ref{thm:main1}, well-known results about the convergence of the Galerkin method based on duality arguments (recapped in Lemma \ref{lem:Schatz} below), and results about the $hp$ approximation spaces in \cite[\S5]{MeSa:10} (used in Lemma \ref{lem:eta} below).

\begin{theorem}[Quasioptimality of the $hp$-FEM if $\Csol(k)$ is polynomially bounded]\label{thm:main2}
Let $d=2$ or $3$, and let $k_0>0$.
Let $(V_N)_{N=0}^\infty$ be the piecewise-polynomial approximation spaces described in \cite[\S5]{MeSa:10} (where, in particular, the triangulations are quasi-uniform), and let $u_N$ be the Galerkin solution defined by \eqref{eq:FEM}.

If $\Csol(k)$ is polynomially bounded (in the sense of Definition \ref{def:polybound}) for $k\in K\subset [k_0,\infty)$ then there exist $C_1,C_2>0$, depending on $\MA, n, R$, and $d$, and $k_0$, but independent of $k$, $h$, and $p$, such that if \eqref{eq:threshold} holds, then, for all $k\in K$, the Galerkin solution exists, is unique, and satisfies the quasi-optimal error bound 
\beq\label{eq:qo}
\N{u-u_N}_{H^1_k(B_R)}\leq \Cqo \min_{v_N\in V_N} \N{u-v_N}_{H^1_k(B_R)},
\eeq
with 
\beq\label{eq:Cqo}
\Cqo:= \frac{2\big(\max\{ A_{\max}, n_{\max} \} + \CDTN_1\big)}{A_{\min}}
\eeq
\end{theorem}

Combining Theorem \ref{thm:main2} with the results on $\Csol(k)$ recapped in Theorem \ref{thm:polybound}, we obtain the following specific examples of coefficients $\MA$ and $n$ when quasioptimality holds.

\begin{corollary}[Quasioptimality under specific conditions on $\MA$ and $n$]\label{cor:qo}
Let $d=2$ or $3$, and let $k_0>0$.

(i) If $\MA$ and $n$ are nontrapping, then there exist $C_1,C_2>0$, depending on $\MA, n, R$, and $d$, and $k_0$, but independent of $k$, $h$, and $p$, such that if \eqref{eq:threshold} holds then, for all $k\geq k_0$, the Galerkin solution exists, is unique, and satisfies the quasi-optimal error bound \eqref{eq:qo} with $\Cqo$ given by \eqref{eq:Cqo}.

(ii) If $\MA$ is $C^\infty$ and $n=1$ then, given $\delta>0$, there exist a set $J$ with $|J|\leq \delta$ and constants $\widetilde{C}_1,\widetilde{C}_2>0$, with all three depending on $\MA, n, R$, $d$, and $k_0$, but independent of $k$, and $\widetilde{C}_2$ additionally depending on $\delta$ and $k_0$ such that, for all $k\in [k_0,\infty)\setminus J$, if \eqref{eq:threshold} holds (with $C_1, C_2$ replaced by $\widetilde{C}_1, \widetilde{C}_2$) then the Galerkin solution exists, is unique, and satisfies \eqref{eq:qo} with $\Cqo$ given by \eqref{eq:Cqo}.
\end{corollary}

For the plane-wave scattering problem (i.e.~for $F(v)$ given by \eqref{eq:planewavedata}), the regularity result 
\beq\label{eq:Cosc}
|u|_{H^2(B_R)} \leq \Cosc k \N{u}_{H^1_k(B_R)}
\eeq
was recently proved in \cite[Theorem 9.1 and Remark 9.10]{LaSpWu:19}, where $\Cosc$ depends on $\MA,n,d,$ and $R$, but is independent of $k$. 
The polynomial approximation bounds in  \cite[\S B]{MeSa:10} imply that, for the sequence of approximation spaces $(V_N)_{N=0}^\infty$ described in \cite[\S5]{MeSa:10}, 
\beq\label{eq:MSbae}
\min_{v_N\in V_N} \N{u-v_N}_{H^1_k(B_R)}\leq C_6 \frac{h}{p} \left( 1 + \frac{kh}{p}\right) |u|_{H^2(B_R)}
\eeq
where $C_6$ only depends on the constants in \cite[Assumption 5.2]{MeSa:10} (which depend on the element maps from the reference element).
Using \eqref{eq:MSbae} and \eqref{eq:Cosc} to bound the right-hand side of \eqref{eq:qo}, we obtain the following bound on the relative error of the Galerkin solution.

\begin{corollary}[Bound on the relative error of the Galerkin solution]
\label{cor:1}
Let the assumptions of Theorem \ref{thm:main2} hold and, furthermore, let $F(v)$ be given by \eqref{eq:planewavedata} (so that $u$ is the solution of the plane-wave scattering problem). 
If $\Csol(k)$ is polynomially bounded (in the sense of Definition \ref{def:polybound}) for $k\in K\subset [k_0,\infty)$, then there exists $C_6>0$, independent of $k$, $h$, and $p$, such that if \eqref{eq:threshold} holds, then, for all $k\in K$,
\beq\label{eq:rel_error}
\frac{
\N{u-u_N}_{H^1_k(B_R)}
}{
\N{u}_{H^1_k(B_R)}
}
\leq \Cqo C_6 \Cosc C_1 \big(1+ C_1\big),
\eeq
with $\Cqo$ given by \eqref{eq:Cqo}; i.e.~the relative error can be made arbitrarily small by making $C_1$ smaller. 
\end{corollary} 

\bre[Theorem \ref{thm:main1} is valid for solutions of a much larger class of PDEs]\label{rem:generalise}
Inspecting the proof of Theorem \ref{thm:main1} below, we see that the conclusion, i.e.~the decomposition 
$u=\uhigh+ \ulow$ with $\uhigh$ and $\ulow$ satisfying the bounds \eqref{eq:uhigh} and \eqref{eq:ulow} respectively, holds under much weaker assumptions. Indeed, the conclusion still holds under the following three assumptions only.

(i) $P_k$ is a family of properly-supported second-order pseudo-differential operators, with principal symbol $p_k(x,\zeta)$,

(ii) $p_k(x,\zeta)$ is coercive at infinity in the sense that
\beq\label{eq:generalise1}
\liminf_{|\xi| \rightarrow \infty, \, x\in \mathbb R^d} \big \langle k\xi \rangle^{-2} p_k(x, k \xi) \geq c > 0,
\eeq
where $c>0$ does not depend on $k$, and 

(iii) the solution to $P_k u = -f$, posed in $\Rea^d$ with $\supp\, f \subset B_R$ and $f\in L^2(B_R)$, satisfies the bound 
$$
\N{u }_{L^2(B_{R+2})} \leq C k^M \Vert f \Vert_{L^2(B_R)},
$$
with $C$ and $M$ independent of $k$, $u$, and $f$. (In fact, the $2$ in the $R+2$ on the left-hand side of the bound can be replaced by any number $>0$.)

In particular, no assumption is made about lower-order terms of $P_k$, or the behaviour of $u$ at infinity (such as a radiation condition). \ere

\section{Recap of relevant results about semiclassical pseudodifferential operators}\label{sec:recap}

The proof of Theorem \ref{thm:main1} relies on standard results about semiclassical pseudodifferential operators. 
We review these here, with our default references being  \cite{Zw:12} and \cite[Appendix E]{DyZw:19}. 
Homogeneous -- as opposed to semiclassical --
versions of the results in this section can be found in, e.g.,
\cite[Chapter 7]{Ta:96}, \cite[Chapter 7]{SaVa:02}, \cite[Chapter
6]{HsWe:08}.\footnote{The counterpart of ``semiclassical''
  involving differential/pseudodifferential operators without a small
  parameter is
  usually called ``homogeneous'' (owing to the homogeneity of the
  principal symbol) rather than ``classical.''  ``Classical'' 
  describes the behaviour in 
either calculus in the small-$\hsc$ or high-frequency limit
respectively, where commutators of operators become Poisson brackets
of symbols, hence classical particle dynamics replaces wave motion.}

While the use of homogeneous pseudodifferential operators in numerical analysis is well established, see, e.g., \cite{SaVa:02},  \cite{HsWe:08}, 
there has been less use of semiclassical pseudodifferential operators.
However, these are ideally-suited for studying the high-frequency behaviour of Helmholtz solutions. Indeed, semiclassical pseudodifferential operators are just pseudodifferential operators with a large/small parameter, 
and behaviour with respect to this parameter is then explicitly kept track of in the associated calculus.

\paragraph{The semiclassical parameter $\hsc= k^{-1}$.}
Instead of working with the parameter $k$ and being interested in the
large-$k$ limit, the semiclassical literature usually works with a
parameter $h:= k^{-1}$ and is interested in the small-$h$ limit. So
that we can easily recall results from this literature, we also work
with the small parameter $k^{-1}$, but to avoid a notational clash
with the meshwidth of the FEM, we let $\hsc:= k^{-1}$ (the notation
$\hsc$ comes from the fact that the semiclassical parameter is related
to Planck's constant, which is written as $2\pi \hsc$; see, e.g., \cite[\S1.2]{Zw:12}, \cite[Page 82]{DyZw:19}, \cite[Chapter 1]{Ma:02a}).
 In this notation, the Helmholtz equation $\nabla\cdot (\MA \nabla u ) + k^2 nu = -f$ becomes 
\beq\label{eq:Helmholtzhsc}
P_{\hsc} u = \hsc^2 f, \quad\text{ where } \quad P_{\hsc} := - \hsc^2 \nabla\cdot(\MA \nabla \cdot) -  n.
\eeq

While some results in semiclassical analysis are valid in the limit $\hsc$ small, the results we recap in this section are valid for all $0<\hsc\leq \hsc_0$ with $\hsc_0<\infty$ arbitrary.

\paragraph{The semiclassical Fourier transform $\mathcal F_{\hsc}$.} 
The semiclassical Fourier transform is defined for $\hsc>0$ by
$$
\mathcal F_{\hsc}\phi(\xi) := \int_{\mathbb R^d} \exp\big( -\ri x \cdot \xi/\hsc\big)
\phi(x) \, \rd x,
$$
and its inverse by
\beq\label{eq:SCFTinverse}
\mathcal F^{-1}_{\hsc}\psi(x) := (2\pi \hsc)^{-d} \int_{\mathbb R^d} \exp\big( \ri x \cdot \xi/\hsc\big)
 \psi(\xi)\, \rd \xi;
\eeq
see \cite[\S3.3]{Zw:12}.
Then 
\beq\label{eq:FTderivative}
\mathcal{F}_\hsc \Big( \big( -\ri \hsc \partial\big)^\alpha \phi\Big) = \xi^\alpha \,\mathcal{F}_\hsc \phi
\eeq
and 
\beq\label{eq:Plancherel}
\N{\phi}_{L^2(\Rea^d)} = \frac{1}{(2\pi \hsc)^{d/2}}\N{\mathcal{F}_\hsc \phi}_{L^2(\Rea^d)}.
\eeq

\paragraph{Semiclassical Sobolev spaces.} 
In the same way that it is convenient to work with the weighted $H^1$
norm \eqref{eq:1knorm} when studying the Helmholtz equation with
parameter $k$, it is convenient to use norms weighted with $\hsc$ when
studying \eqref{eq:Helmholtzhsc}. Therefore on the space
\beqs
H_\hsc ^ s (\Rea^d):= \Big\{ u\in L^2(\mathbb R^d), \; \langle \xi \rangle^s 
\mathcal F_\hsc u \in  L^2(\mathbb R^d) \Big\}, \quad\text{ where }\langle \xi \rangle := (1+|\xi|^2)^{1/2}, \quad s\in \Rea,
\eeqs
we use the norm
\beq\label{eq:Hhnorm}
\Vert u \Vert_{H_\hsc^s(\Rea^d)} ^2 := (2\pi \hsc)^{-d} \int_{\Rea^d} \langle \xi \rangle^{2s}
 |\mathcal F_\hsc u(\xi)|^2 \, \rd \xi;
\eeq
see \cite[\S8.3]{Zw:12}, \cite[\S E.1.8]{DyZw:19}. 
We abbreviate $H_\hsc ^ s (\Rea^d)$ to $H_\hsc ^ s$ and $L^2(\Rea^d)$ to $L^2$.

We record for later the fact that, by \eqref{eq:FTderivative} and \eqref{eq:Plancherel}, for multiindices $\alpha$,
\begin{align}
\hsc^{|\alpha|} \N{\partial^\alpha \phi}_{L^2} &= 
\N{ \big(-\ri \hsc \partial\big)^\alpha \phi}_{L^2} = 
\frac{1}{(2\pi \hsc)^{d/2}} \N{\xi^\alpha \, \mathcal{F}_\hsc \phi}_{L^2} \leq 
\frac{1}{(2\pi \hsc)^{d/2}} \big\|\langle\xi\rangle^{|\alpha|} \, \mathcal{F}_\hsc \phi\big\|_{L^2} = 
 \N{\phi}_{H^{|\alpha|}_\hsc}.
\label{eq:H2FT}
\end{align}

\paragraph{Phase space.}
The set of all possible positions $x$ and momenta (i.e.~Fourier variables) $\xi$ is denoted by $T^*\Rea^d$; this is known informally as ``phase space". Strictly, $T^*\Rea^d :=\Rea^d \times (\Rea^d)^*$, but 
for our purposes, we can consider $T^*\Rea^d$ as $\{(x,\xi) : \bx\in \Rea^d, \xi\in\Rea^d\}$.

To deal with the behavior of
functions on phase space uniformly near $\xi=\infty$ (so-called \emph{fiber infinity}), we consider the \emph{radial
  compactification} in the $\xi$ variable of $T^*\Rea^d$. This is defined by
$$
\overline{T} ^* \mathbb R^d:= \mathbb R^d \times B^d,
$$
where $B^d$ denotes the closed unit ball, considered as the closure of the
image of $\mathbb R^d$ under the radial compactification map 
$$\RC: \xi \mapsto \xi/(1+\langle \xi
\rangle);$$
see \cite[\S E.1.3]{DyZw:19}.
Near the boundary of the
ball, $\lvert \xi\rvert^{-1}\circ \RC^{-1}$ is a smooth function, vanishing to
first order at the boundary, with $(\lvert \xi\rvert^{-1}\circ \RC^{-1}, \widehat\xi\circ\RC^{-1})$
thus giving local coordinates on the ball near its boundary.  The boundary of the
ball should be considered as a sphere at infinity consisting of all
possible \emph{directions} of the momentum variable.  More generally, we denote $\overline{T}^* X := X \times B^d$ for $X \subset \mathbb R^d$, and where
appropriate (e.g., in dealing with finite values of $\xi$ only), we abuse notation by dropping the composition with $\RC$ from our
notation and simply identifying $\mathbb R^d$ with the interior of $B^d$.

\paragraph{Symbols, quantisation, and semiclassical pseudodifferential operators.} 

A symbol is a function on $T^*\Rea^d$ that is also allowed to depend on $\hsc$, and thus can be considered as an $\hsc$-dependent family of functions.
Such a family $a=(a_\hsc)_{0<\hsc\leq\hsc_0}$, with $a_\hsc \in C^\infty({T^*\mathbb R^d})$, 
is a \emph{symbol
of order $m$}, written as $a\in S^m(\Rea^d)$, if for any multiindices $\alpha, \beta$
\beq\label{eq:Sm}
| \partial_x^\alpha \partial^\beta_\xi a(x,\xi) | \leq C_{\alpha, \beta} \langle \xi \rangle^{m -|\beta|} \quad\tfa (x,\xi) \in T^* \Rea^d \text{ and for all } 0<\hsc\leq \hsc_0,
\eeq
where $C_{\alpha, \beta}$ does not depend on $\hsc$, $x$, or $\xi$; see \cite[p.~207]{Zw:12}, \cite[\S E.1.2]{DyZw:19}. 
In this paper, we only consider these symbol classes on $\Rea^d$, and so we abbreviate $S^m(\Rea^d)$ to $S^m$.

For $a \in S^m$, we define the \emph{semiclassical quantisation} of $a$, $\operatorname{Op}_{\hsc}(a): \mathscr{S}(\Rea^d)\rightarrow \mathscr{S}(\Rea^d)$, by  
\beq \label{Oph}
\big(\operatorname{Op}_{\hsc}(a) v\big)(x) := (2\pi \hsc)^{-d} \int_{\Rea^d} \int_{\Rea^d} 
\exp\big(\ri (x-y)\cdot\xi/\hsc\big)\,
a(x,\xi) v(y) \,\rd y  \rd \xi
\eeq
for $v\in \mathscr{S}(\Rea^d)$; \cite[\S4.1]{Zw:12} \cite[Page 543]{DyZw:19}. The integral in \eqref{Oph} need not converge, and can be understood \emph{either} as an oscillatory integral in the sense of \cite[\S3.6]{Zw:12}, \cite[\S7.8]{Ho:83}, \emph{or} as an iterated integral, with the $y$ integration performed first; see \cite[Page 543]{DyZw:19}.

Conversely, if $A$ can be written in the form above, i.\,e.\ $A = \operatorname{Op}_{\hsc}(a)$ with $a\in S^m$, we say that $A$ is a \emph{semiclassical pseudo-differential operator of order $m$} and
we write $A \in \Psi_{\hsc}^m$. We use the notation $a \in \hsc^l S^m$  if $\hsc^{-l} a \in S^m$; similarly 
$A \in \hsc^l \Psi_\hsc^m$ if 
$\hsc^{-l}A \in \Psi_\hsc^m$.

\begin{theorem}\mythmname{Composition and mapping properties of
semiclassical pseudo-differential operators \cite[Theorem 8.10]{Zw:12}, \cite[Proposition E.17 and Proposition E.19]{DyZw:19}}\label{thm:basicP} 
If $A\in \Psi_{\hsc}^{m_1}$ and $B  \in \Psi_{\hsc}^{m_2}$, then
\begin{itemize}
\item[(i)]  $AB \in \Psi_{\hsc}^{m_1+m_2}$,
\item[(ii)]  $[A,B]:= AB-BA \in \hsc\Psi_{\hsc}^{m_1+m_2-1}$,
\item[(iii)]  For any $s \in \mathbb R$, $A$ is bounded uniformly in $\hsc$ as an operator from $H_\hsc^s$ to $H_\hsc^{s-m_1}$.
\end {itemize}
\end{theorem}

\paragraph{Residual class.} 
We say that $A =O(\hsc^\infty)_{\Psi^{-\infty}}$ if, for any $s>0$ and $N\geq 1$, there exists $C_{s,N}>0$ so that
\beq \label{eq:residual}
\Vert A \Vert_{H_\hsc^{-s} \rightarrow H_\hsc^{s}} \leq C_{N,s} \hsc^N;
\eeq
i.e.~$A\in \Psi_\hsc^{-\infty}$ and furthermore all of its operator norms are bounded by any algebraic power of $\hsc$.

\paragraph{Principal symbol $\sigma_{\hsc}$.}
Let the quotient space $ S^m/\hsc S^{m-1}$ be defined by identifying elements 
of  $S^m$ that differ only by an element of $\hsc S^{m-1}$. 
For any $m$, there is a linear, surjective map
$$
\sigma^m_{\hsc}:\Psi_\hsc ^m \to S^m/\hsc S^{m-1},
$$
called the \emph{principal symbol map}, 
such that, for $a\in S^m$,
\beq\label{eq:symbolone}
\sigma_\hsc^m\big(\Op_\hsc(a)\big) = a \quad\text{ mod } \hsc S^{m-1};
\eeq
see \cite[Page 213]{Zw:12}, \cite[Proposition E.14]{DyZw:19} (observe that \eqref{eq:symbolone} implies that 
$\operatorname{ker}(\sigma^m_{\hsc}) = \hsc\Psi_\hsc ^{m-1}$).

When applying the map $\sigma^m_{\hsc}$ to 
elements of $\Psi^m_\hsc$, we denote it by $\sigma_{\hsc}$ (i.e.~we omit the $m$ dependence) and we use $\sigma_{\hsc}(A)$ to denote one of the representatives
in $S^m$ (with the results we use then independent of the choice of representative).
Key properties of the principal symbol that we use below are that 
\beq \label{eq:multsymb}
\sigma_{\hsc}(AB)=\sigma_{\hsc}(A)\sigma_{\hsc}(B),
\eeq
\beq\label{eq:symbolPDE}
\sigma_{\hsc}(P_{\hsc})=  \langle \MA \xi, \xi\rangle - n,
\eeq
where $\langle\cdot,\cdot\rangle$ denotes the $\ell^2$ inner product on $\Rea^d$. 
%
The property 
\eqref{eq:multsymb} is proved in \cite[Proposition E.17]{DyZw:19}, \eqref{eq:symbolPDE} follows from \eqref{eq:symbolone} since $P_\hsc = \Op_\hsc \big(\langle \MA \xi, \xi\rangle - n  - \ri \hsc \xi_\ell \partial_j A_{j\ell}\big)$ (where we sum over the indices $j$ and $\ell$).

\paragraph{Operator wavefront set $\operatorname{WF}_{\hsc}$.}  
We say that $(x_0,\xi_0) \in {\overline{T}}^*\mathbb{R}^d$ is \emph{not} in the \emph{semiclassical operator wavefront set} of $A = \operatorname{Op}_{\hsc}(a) \in \Psi_{\hsc}^m$, denoted by $\operatorname{WF}_{\hsc} A$, if there exists a neighbourhood $U$ of $(x_0,\xi_0)$ such that for all multiindices $\alpha, \beta$ and all $N\geq 1$ there exists $C_{\alpha,\beta,U,N}>0$ (independent of $\hsc$) so that, for all $0<\hsc\leq \hsc_0$,
\beq\label{eq:microsupport}
|\partial_x^\alpha \partial_\xi^\beta  a(x,\xi)| \leq C_{\alpha,
  \beta, U, N} \hsc^N \langle \xi \rangle^{-N}\quad\tfa (x,\RC(\xi))\in U;
\eeq
i.e.~outside its semiclassical operator wavefront set an operator
vanishes faster than any algebraic power of both $\hsc$ and
  $\langle \xi \rangle^{-1}$; see \cite[Page 194]{Zw:12}, \cite[Definition E.27]{DyZw:19}.
Three properties of the semiclassical operator wavefront set that we use below are 
\beq \label{eq:WFprod}
\operatorname{WF}_{\hsc} (AB) \subset \operatorname{WF}_{\hsc} A \cap \operatorname{WF}_{\hsc}B
\eeq
(see \cite[\S 8.4]{Zw:12}, \cite[E.2.5]{DyZw:19}),
\beq\label{eq:support}
\operatorname{WF}_{\hsc}\big( \Op_\hsc(a)\big) \subset \supp \,a
\eeq
(since $(\supp \, a)^c \subset (\operatorname{WF}_{\hsc}( \Op_\hsc(a)))^c$ by \eqref{eq:microsupport}), and 
\beq\label{eq:ResWF}
\operatorname{WF}_\hbar A = \emptyset \iff A = O(\hbar^\infty)_{\Psi^{-\infty}}
\eeq
(see \cite[E.2.2]{DyZw:19}).

\paragraph{Compactly-supported operators.}
We say that  $A$ is \emph{compactly supported} if its Schwartz kernel
is compactly supported in some set $K \Subset \Rea^d \times \Rea^d,$
for all $0<\hsc\leq \hsc_0.$
We recall that if $\mathcal{D}(\Rea^d):= C^{\infty}_{\rm comp}(\Rea^d)$ (i.e.~the set of test functions) and $\mathcal{D}'(\Rea^d)$ denote the set of linear functionals on $\mathcal{D}(\Rea^d)$ (i.e.~the set of distributions),
given a bounded, sequentially-continuous operator $A : \mathcal D \rightarrow\mathcal D'$ there exists a \emph{Schwartz kernel} $\mathcal K_A\in\mathcal{D}'(\Rea^d\times\Rea^d) $ such that
$$
A v (x) = \int_{\Rea^d} \mathcal K_A(x,y)v(y) \,\rd y,
$$
in the sense of distributions; see, e.g., \cite[Theorem 5.2.1]{Ho:83}, \cite[\S A.7]{DyZw:19}. We use below the facts that 
\bit
\item $A$ is compactly supported iff there exist $\chi_1,\chi_2\in \mathcal{D}$ such that $A=\chi_1 A \chi_2$, thus
\item if $\chi_1, \chi_2 \in \mathcal{D}$ are compactly supported functions, then $\chi_1 A \chi_2$ is compactly supported, and
\item if  $P$ is a differential operator and $\chi \in \mathcal{D}$, then both $\chi P$ and $P\chi $ are compactly supported.
\eit

\paragraph{Ellipticity.} 
We  say that $B\in \Psi_\hsc^m$ is \emph{elliptic} on $X \subset \overline{T}^*\mathbb R^d$ if there exists $c>0$, independent of $\hsc$, such that 
\beq \label{eq:ellip}
\langle \xi \rangle^{-m} \big|\sigma_\hsc(B)(x,\xi)\big| \geq c, \quad \tfa (x,\RC(\xi))\in X  \text{ and for all } 0<\hsc\leq \hsc_0.
\eeq

A key feature of elliptic operators is that they are microlocally invertible; this is reflected in the following result.

\begin{proposition}\mythmname{Elliptic parametrix \cite[Proposition
    E.32]{DyZw:19}} \footnote{\label{foot:proper}We highlight that working in $\mathbb R^d$ (as opposed to on a general manifold defined by coordinate charts) allows us to remove the proper-support assumption appearing in \cite[Proposition E.32, Theorem E.33]{DyZw:19}.} 
 \label{prop:para} 
Let $A \in \Psi_\hsc^{m}$ and $B \in \Psi_\hsc^{\ell}$ be such 
that $B$ is elliptic on $\operatorname{WF_\hsc}(A)$.
Then there exist $Q, Q' \in \Psi_\hsc^{m-\ell}$ such that
$$
A = BQ + O(\hsc^\infty)_{\Psi^{-\infty}} = Q'B + O(\hsc^\infty)_{\Psi^{-\infty}}.
$$ 
\end{proposition}

\begin{theorem}\mythmname{Elliptic estimate \cite[Theorem E.33]{DyZw:19}} 
$^{\ref{foot:proper}}$ \label{thm:elliptic}
Let $A \in \Psi_\hsc^{m_1}$, $B_1 \in \Psi_\hsc^{m_2}$, and $P\in \Psi_\hsc^\ell$ be so that $B_1 P$ is elliptic on $\WFh (A)$.

(i) Given $s, N>0,$ and $M>0$, if $v\in \mathcal{D}'$ and $B_1 P v\in H^{s-m_2-\ell}$ then $A v \in H^{s-m_1}$ and there exists $C_s>0$, $C_{N,M,s}>0$ (independent of $v$ and $\hsc$)  such that 
\beq\label{eq:elliptic_estimate}
\N{A v}_{H^{s-m_1}_\hsc} \leq C_s \N{B_1 P v}_{H^{s-m_2-\ell}_\hsc} + C_{N,M,s} \, \hsc^M \N{v}_{H^{-N}_\hsc}.
\eeq

(ii) If, in addition, $A$ and $B_1 P$ are compactly supported, then there exists $\widetilde{\chi}\in C^\infty_{\rm comp}$ so that
\beq\label{eq:elliptic_estimate_comp}
\N{A v}_{H^{s-m_1}_\hsc} \leq C_s \N{B_1 P v}_{H^{s-m_2-\ell}_\hsc} + C_{N,M,s} \, \hsc^M \N{\widetilde{\chi} v}_{H^{-N}_\hsc}.
\eeq
\end{theorem}

Part (i) of Theorem \ref{thm:elliptic} is proved by using Proposition \ref{prop:para} with $B= B_1 P\in \Psi_\hsc^{m_2+\ell}$, applying the resulting operator equation to $v$, and taking norms. The operator $Q' \in \Psi_\hsc^{m_1-m_2-\ell}$ and the constant $C_s$ is then $\|Q'\|_{H^{s-m_2-\ell}_\hsc \rightarrow H^{s-m_1}_\hsc}$. The proof of Part (ii) is similar, using that, since $A$ and $B_1P$ are both compactly supported, there exists $\widetilde{\chi}\in C^\infty_{\rm comp}$ such that $(A-B_1P)v = (A-B_1P)\widetilde{\chi}v$.

\section{Proof of Theorem \ref{thm:main1}}\label{sec:proof_PDE}

In the notation introduced in \S\ref{sec:recap}, Theorem \ref{thm:main1} becomes the following.

\begin{theorem} \label{th:mainsemi}
Let $\MA$ and $n$ satisfy Assumption \ref{ass:1} and let $R>0$ be such that $\supp(\MI-\MA)\cup\supp(1-n)\Subset B_R$. Given $f\in L^2(B_R)$, let $u$ satisfy $P_\hsc u=\hsc^2 f$ in $\Rea^d$ and the Sommerfeld radiation condition \eqref{eq:src}.
Assume that, given $k_0>0$, $\Csol(k)$ is polynomially bounded (in the sense of Definition \ref{def:polybound}) for $k\in K \subset [k_0,\infty)$.
Given $k_0>0$, let $\hsc_0:= k_0^{-1}$, and let $H:=\{ k^{-1} : k \in K\} \subset (0, \hsc_0]$.

Then there exist $C_3, C_4, C_5>0$ such that
\beqs
u|_{B_R}= \uhigh+ \ulow
\eeqs
where $\uhigh \in H^2_\hsc(B_R)$ with
\beq\label{eq:uhigh_bound_proof}
\N{\partial^\alpha \uhigh}_{L^2(B_R)} \leq C_3 \hsc^{2-|\alpha|} \N{f}_{L^2(B_R)}
 \quad\text{ for all $|\alpha|\leq 2$ and for all $\hsc\in H\subset (0,\hsc_0]$},
\eeq
and $\ulow \in C^\infty(B_R)$ with
\beq\label{eq:ulow_bound_proof}
\N{\partial^\beta \ulow}_{L^2(B_R)} \leq \Csol\big(\hsc^{-1}; R+2\big) \,\, C_4\, \left(\frac{\hsc}{C_5}\right)^{1-|\beta|} \N{f}_{L^2(B_R)} \,\,\text{ for all $\beta$ and for all $
\hsc\in H\subset (0,\hsc_0]$},
\eeq
where $C_3, C_4,$ and $C_5$ depend on $\MA, n, d$, and $\hsc_0$, but are independent of $\hsc$, $f$, $\alpha$, and $\beta$.
\end{theorem}

\subsection{Step 0: Restatement of bounds on the solution operator in semiclassical notation}

The definition of $\Csol$ (Definition \ref{def:Csol}) implies that, in semiclassical notation,
\beq\label{eq:Csolh}
\N{u}_{H^1_\hsc(B_R)} \leq \hsc\, \Csol(\hsc^{-1})\N{f}_{L^2(B_R)} \quad\tfa \hsc>0.
\eeq
It is convenient to record here in semiclassical notation the bound on the solution operator when $\Csol$ is polynomially bounded.

\ble[Polynomial boundedness rewritten in terms of $\hsc$]
Given $f \in L^2_{\comp}(\Rea^d)$, let $u\in H^1_{\rm loc}(\Rea^d)$ be the solution to 
\beqs
P_{\hsc} u= \hsc^2 f
\eeqs
satisfying the Sommerfeld radiation condition \eqref{eq:src} (with $k = \hsc^{-1}$).

If $\Csol(k)$ is polynomially bounded for $k\in K \subset [k_0,\infty)$ (in the sense of Definition \ref{def:polybound}), then there exists $M>0$ (independent of $\hsc$) such that, given $\chi \in C^\infty_{\rm comp}(\Rea^d)$, there exists $C>0$ (independent of $\hsc$ but dependent on $\chi$) such that
\beq\label{eq:res}
\N{\chi u}_{L^2} \leq C\hsc^{1-M} \N{f}_{L^2} \quad\tfor \hsc \in H \subset (0,\hsc_0], 
\eeq
where $\hsc_0:= k_0^{-1}$ and $H:=\{ k^{-1} : k \in K\}$.
\ele
The bound \eqref{eq:res} also holds with $\N{\chi u}_{L^2}$ replaced by $\N{\chi u}_{H^1_\hsc}$, but we only need it in the form \eqref{eq:res} for what follows.

\subsection{Step 1: The definitions of $\ulow$ and $\uhigh$.}

\paragraph{The cut-off functions $\chi$ and $\chi_\mu$.}
Let
$\chi \in C^\infty_{\rm comp} (\mathbb R^d; [0,1])$ be such that 
\beq\label{eq:chi}
\chi = 
\begin{cases}
1 &\text{ in } B_1 \\
0 &\text{ outside } B_2.
\end{cases}
\eeq
For $\mu>0$, let
\beq\label{eq:chilambda}
\chi_\mu(\cdot) := \chi\left(\frac{\cdot }{\mu}\right).
\eeq
We define $\mu_0= \mu_0(\MA,n)$ by 
\beq\label{eq:mu0}
\mu_0(\MA,n):= \left( 1+\frac{2n_{\max}}{A_{\min}}\right).
\eeq
The reason for this definition is that it implies that
\beq\label{eq:lambda0}
\text{ if } \quad |\xi|^2\geq \mu_0 \quad\text{ then }  \quad\langle \xi \rangle^{-2} \sigma_{\hsc}(P) \geq \frac{A_{\min}}{2}>0.
\eeq
Indeed, by \eqref{eq:symbolPDE},
\beqs
\langle \xi \rangle^{-2} \sigma_{\hsc}(P) \geq \frac{A_{\min} |\xi|^2 -n_{\max}}{1+|\xi|^2} = \frac{A_{\min}}{2}
+\left(\frac{A_{\min}}{2}\right) 
\left(\frac{
|\xi|^2 - 1 - 2 n_{\max}/A_{\min}
}{
1+ |\xi|^2
}\right),
\eeqs
and \eqref{eq:lambda0} follows. The importance of the property \eqref{eq:lambda0} is explained at the end of this subsection.

\paragraph{The frequency cut-offs $\Pilow$ and $\Pihigh$.}

We define $\Pilow$ and $\Pihigh$, the projections on low and high frequencies respectively,
by \eqref{eq:Pilow_intro} and \eqref{eq:Pihi_intro}. The definition of the quantisation $\operatorname{Op}_\hsc$ \eqref{Oph} and the change of variable $\zeta = \xi/\hsc$ imply that 
\beq \label{Low}
\Pilow= \operatorname{Op}_{\hsc}\big(\chi_\mu(|\xi|^2) \big)
\eeq
and
\beq \label{High}
\Pihigh=I-\Pilow. 
\eeq
These definitions and the definition of $\Psi^m_\hsc(\Rea^d)$ in \S\ref{sec:recap} imply that $\Pilow \in \Psih^{-\infty}(\RR^d)$ and $\Pihigh\in \Psih^{0}(\RR^d).$

\paragraph{The locations of the wavefront sets of the frequency cut-offs, and the regions where their symbols equal one.}

In Figure \ref{fig:1} we show, as functions of $|\xi|^2$,  the locations of 
$\operatorname{WF}_\hsc (\Pihigh) $ and $\operatorname{WF}_\hsc (\Pilow)$, and the regions where $\sigma_\hsc(\Pihigh),$ and $\sigma_\hsc(\Pilow)$ equal one.
These locations/regions are obtained using \eqref{eq:support} and \eqref{eq:symbolone} respectively.
For example, since $1 -\chi_\mu(|\xi|^2) = 1$ for $|\xi|^2\geq 2\mu$ and $=0$ for $|\xi|^2 \leq \mu$, \eqref{eq:symbolone} and \eqref{eq:support} imply that 
\beq\label{eq:symbolWFkey}
\sigma_\hsc(\Pihigh) = 1 \text{ on } \big\{ \xi \, :\, |\xi|^2 \geq 2\mu\big\}\quad\tand\quad 
\WFh(\Pihigh) \subset \big\{ \xi \, :\, |\xi|^2 \geq \mu\big\}.
\eeq
We also record the following key consequence of the results summarised in Figure \ref{fig:1}.
\ble\label{lem:Pelliptic}
If $\mu \geq \mu_0$, then $P_\hsc$ is elliptic on $\WFh(\Pihigh)$.  
\ele

This property is central to our proof of the bound \eqref{eq:uhigh_bound_proof} on $\uhigh$, i.e., the high-frequency component.
It is a consequence of \eqref{eq:lambda0}, and the reason why we choose $\mu_0$ as in \eqref{eq:mu0} is for this ellipticity result to hold. 

\begin{figure}
\centering{
\includegraphics[width=0.8\textwidth]{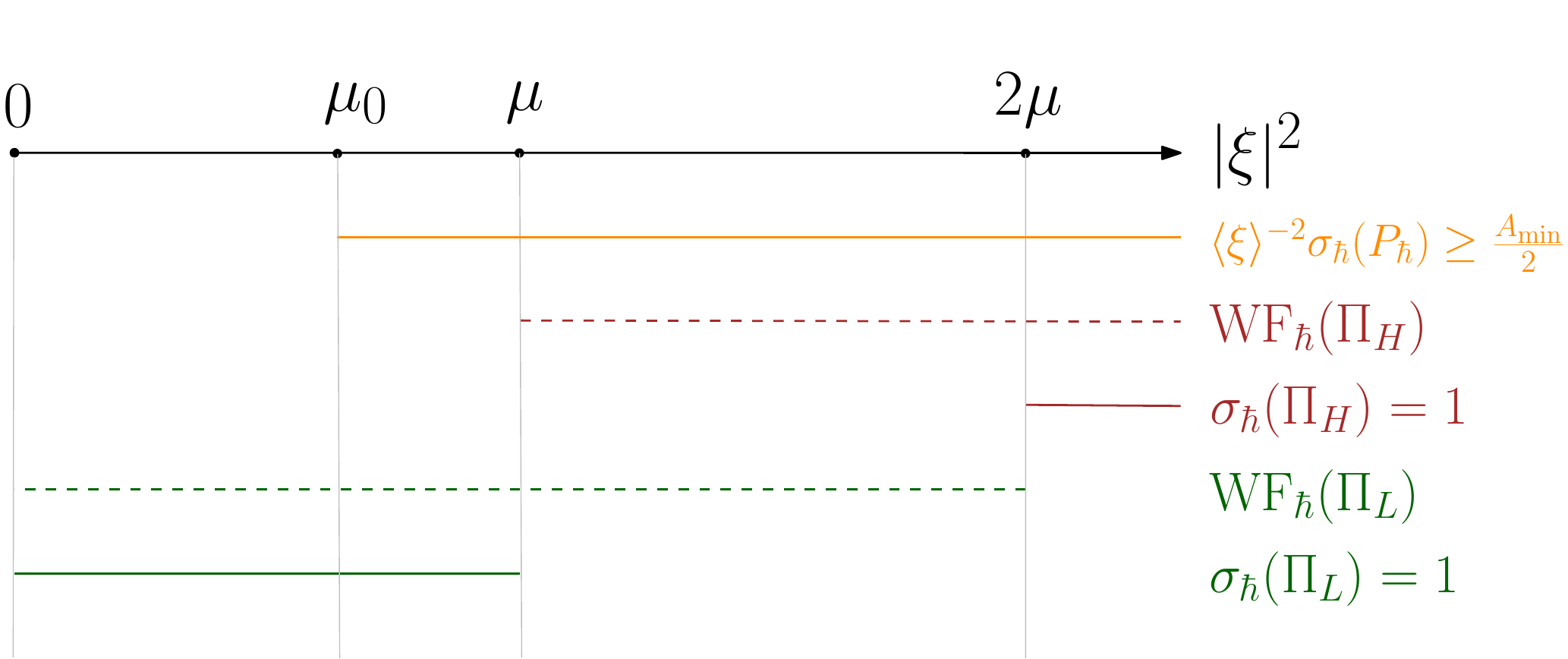}
}
\caption{The locations of $\operatorname{WF}_\hsc (\Pihigh) $ and $\operatorname{WF}_\hsc (\Pilow)$, the regions where the principal symbols of $\Pihigh$ and $\Pilow$ equal one, and the region where $P_\hbar$ is elliptic.}\label{fig:1}
\end{figure}

\paragraph{The definitions of $\ulow$ and $\uhigh$.}

As described in \S\ref{sec:informal}, we choose $\varphi \in C^\infty_{\rm comp}(\Rea^d)$ be equal to one on $B_{R+1}$ and vanish outside $B_{R+2}$. 
We then let 
$$
w := \varphi u
$$
and we define 
$$
u_{\mathcal A} := (\Pi_L w)\big|_{B_R} \quad\tand\quad u_{H^2} := (\Pi_H w)\big|_{B_R}.
$$

\subsection{Step 2: Proof of the bound \eqref{eq:ulow_bound_proof} on $\ulow$ (the low-frequency component)}\label{sec:lowproof}

Since $\Pilow \in \Psi_\hsc^{-\infty}$, Part (iii) of Theorem \ref{thm:basicP}, together with Sobolev embedding, gives $\Pilow w \in C^\infty$.

The definition of $\Pilow$ \eqref{eq:Pilow_intro} and Plancherel's identity \eqref{eq:Plancherel} for the standard (i.e.~non semiclassical) Fourier transform imply that
\beq\label{eq:final1}
  \N{\pa^\beta\big( \Pilow w\big)}_{L^2}=\frac{1}{(2\pi)^{d/2}}\N{(\cdot)^\beta \mathcal{F}\big( \Pilow w\big)(\cdot)}_{L^2} 
=  \frac{1}{(2\pi)^{d/2}}\N{(\cdot)^\beta \chi_\mu\big(\hsc^2 |\cdot|^2\big)\mathcal{F}w(\cdot)}_{L^2}.
  \eeq
The definitions of $\chi$ \eqref{eq:chi} and $\chi_\mu$ \eqref{eq:chilambda} imply that $\chi_\mu(\xi)=0$ for $|\xi|\geq 2 \mu$, so
\beqs
\chi_\mu\big(\hsc^2 |\zeta|^2\big) = 0 \quad\tfor \, |\zeta| \geq \sqrt{2\mu} \, \hsc^{-1}.
\eeqs  
Using this fact, and then (in this order) the fact that $|\chi_\mu|\leq 1$, Plancherel's identity for the standard Fourier transform, the fact that $\varphi=0$ outside $B_{R+2}$, and the definition of $\Csol$ \eqref{eq:Csol}, we find from \eqref{eq:final1} that
\begin{align*}\nonumber
  \N{\pa^\beta\big( \Pilow \varphi u\big)}_{L^2}
  &\leq \frac{(2 \mu)^{|\beta|/2}}{(2\pi)^{d/2}}  \hsc^{-|\beta|} \N{\chi_\mu\big(\hsc^{2} |\cdot|^2\big)\mathcal{F}(\varphi u)(\cdot)}_{L^2}\\
  &\leq \frac{(2 \mu)^{|\beta|/2}}{(2\pi)^{d/2}}  \hsc^{-|\beta|} \N{\mathcal{F}(\varphi u)}_{L^2}\\
    &\leq (2 \mu)^{|\beta|/2}  \hsc^{-|\beta|} \N{\varphi u}_{L^2}\\
        &\leq {(2 \mu)^{|\beta|/2}  \hsc^{-|\beta|} \hsc\,\Csol(\hsc^{-1}; R+2)} \N{f}_{L^2(B_R)}.
  \end{align*}
Since   
\beqs
\N{\partial^\beta \ulow}_{L^2(B_R)} = \N{\partial^\beta (\Pilow w)}_{L^2(B_R)} \leq  \N{\partial^\beta (\Pilow w)}_{L^2},
\eeqs
  the bound \eqref{eq:ulow_bound_proof} then follows with  $C_4:= \sqrt{2\mu}$ and $C_5:= \sqrt{2\mu}$.

\subsection{Step 3: Proof of the bound \eqref{eq:uhigh_bound_proof} on $\uhigh$ (the high-frequency component)}\label{sec:5.4}

By the inequality \eqref{eq:H2FT},
it is sufficient to prove  that 
\beq\label{eq:HF2}
\N{\Pihigh w}_{H^2_\hsc}\leq C_3 \hsc^{2} \N{f}_{L^2(B_R)} \quad\text{ for all $\hsc\in H\subset (0,\hsc_0)$}.
\eeq

It is instructive to first prove \eqref{eq:HF2} under the assumption that $\Csol(k)\lesssim 1$ (which, by Theorem \ref{thm:polybound} is ensured if $\MA$ and $n$ are nontrapping). Indeed, as discussed in \S\ref{sec:informal}, this proof only requires 
that $P_\hsc$ is elliptic on $\WFh(\Pihigh)$; i.e.,~Lemma \ref{lem:Pelliptic}.
Throughout the rest of this section, therefore, we assume that $\mu\geq \mu_0$, so that the result of Lemma \ref{lem:Pelliptic} holds.

\subsubsection{Proof of \eqref{eq:HF2} under the assumption that $\Csol(k)\lesssim 1$}\label{sec:simpleproof}

We seek to apply Part (i) of Theorem \ref{thm:elliptic} with $A= \Pi_H$ (so $m_1=0$), $B_1 =1$ (so $m_2=0$), and $P=P_\hsc$ (so $\ell=2$).
By Lemma \ref{lem:Pelliptic}, $B_1 P$ is elliptic on $\WFh (A)$.
We can therefore apply Theorem \ref{thm:elliptic} and obtain that, given $N, N'>0$,
\beq\label{eq:Friday1}
\big\|\Pi_H w \big\|_{H^2_\hsc} \lesssim \N{ P_\hsc w}_{L^2} + \hsc^{N'} \N{w}_{H^{-N}_\hsc},
\eeq
where the omitted constant in $\lesssim$ depends on $N$ and $N'$.
Since $P_\hsc u = \hsc^2 f$,
\beqs
P_\hsc w = [P_\hsc, \varphi]u+\hsc^2\varphi f,
\eeqs
where $[\cdot,\cdot]$ is the standard commutator defined by $[A_1, A_2]:=A_1 A_2 - A_2 A_1$, 
so that \eqref{eq:Friday1} becomes
\beq\label{eq:Friday1b}
\big\|\Pi_H w \big\|_{H^2_\hsc} \lesssim \N{  [P_\hsc, \varphi]u }_{L^2}+ \hsc^2\N{f}_{L^2} + \hsc^{N'} \N{w}_{H^{-N}_\hsc}.
\eeq
Direct calculation, using the fact that $\supp \,\varphi \subset B_{R+2}$, implies that
\beq\label{eq:Friday2}
\big\| [P_\hsc,\varphi]u\big\|_{L^2} \lesssim \hsc \N{u}_{H^1_\hsc(B_{R+2})},
\eeq
where the omitted constant depends on $\varphi$, and hence on $R$.

Combining \eqref{eq:Friday1b} 
and \eqref{eq:Friday2}, and recalling that $\supp\, \varphi \subset B_{R+2}$, we have
\beqs
\big\|\Pi_H w \big\|_{H^2_\hsc} \lesssim \hsc \N{u}_{H^1_\hsc(B_{R+2})} + {\hsc^2}\N{f}_{L^2(B_R)}+ \hsc^{N'} \N{u }_{H^{-N}_\hsc(B_{R+2})}.
\eeqs
Choosing $N=0$ and $N'=1$, and then using \eqref{eq:Csolh},
we obtain 
\beq\label{eq:Sunday1}
\big\|\Pi_H w \big\|_{H^2_\hsc} \lesssim \hsc^2 \Big( 1 + \Csol(\hsc^{-1})\Big) \N{f}_{L^2(B_R)}.
\eeq
If $\Csol (\hsc^{-1})\lesssim 1$, then this implies \eqref{eq:HF2}. 
However, if $\Csol(\hsc^{-1})\gg 1$ (as occurs when $\Csol$ is polynomially bounded in the sense of Definition \ref{def:polybound} with $M>0$) then \eqref{eq:Sunday1} is a weaker bound than \eqref{eq:HF2}.

\subsubsection{Proof of \eqref{eq:HF2} under the assumption that $\Csol(k)$ is polynomially bounded}\label{sec:betterproof}
 Inspecting the argument in \S\ref{sec:simpleproof}, we see that the assumption that $\Csol(k)\lesssim 1$ is needed to get a good bound on the commutator term $[P_\hsc,\varphi]u$. To remove this commutator term, one idea is to use the elliptic estimate in Part (i) of Theorem \ref{thm:elliptic}, using the fact that $P_\hbar$ is elliptic on $\operatorname{WF}_\hbar( \Pi_H \varphi)$, and apply the estimate with $v:=u$. However, the error term would not be compactly supported and we would be unable to control it using the polynomial bound on the solution operator \eqref{eq:res}. We therefore introduce additional spatial cut-offs on the left of $\Pi_H \varphi$ and $P_\hbar$ to create compactly-supported operators and have a compactly-supported error term thanks to Part (ii) of Theorem \ref{thm:elliptic}.

To this end, let $\varphi_1, \; \varphi_2 \in C^\infty_{\rm comp}(\mathbb R^d)$ be such that $\varphi_1 = 1$ on $\operatorname{supp} \varphi$ and  $\varphi_2 = 1$ on $\operatorname{supp} \varphi_1$; we then write
\beq \label{eq:rev_1}
\Pi_H \varphi u = (1-\varphi_1)\Pi_H \varphi u + \varphi_1 \Pi_H \varphi u.
\eeq
Since $1-\varphi_1 = 0$ on $\operatorname{supp}\varphi$, using (\ref{eq:WFprod}) and (\ref{eq:support}), we obtain that
$$
\operatorname{WF}_\hbar \big((1-\varphi_1)\Pi_H \varphi \big)\subset \overline{T}^*( \operatorname{supp}(1-\varphi_1)) \cap \overline{T}^* (\operatorname{supp} \varphi) = \emptyset.
$$
Hence, by (\ref{eq:ResWF}), $(1-\varphi_1)\Pi_H \varphi = O(\hbar^\infty)_{\Psi^{-\infty}}$, and, by the definition of the residual class (\ref{eq:residual}), for any $N \geq 1$ there exists $C_N >0$ so that
\beq \label{eq:rev_2}
\Vert (1-\varphi_1)\Pi_H \varphi u \Vert_{H^2_\hbar} = \Vert (1-\varphi_1)\Pi_H \varphi \varphi_1 u \Vert_{H^2_\hbar} \leq C_N \hbar^N \Vert \varphi_1 u\Vert_{L^2},
\eeq
were we used the fact that $\varphi_1 = 1$ on $\operatorname{supp}\varphi$ in the first equality.

It therefore remains to control $\varphi_1\Pi_H \varphi u$; to do this, we use the elliptic estimate of Theorem \ref{thm:elliptic}.
\ble\label{lem:rev_ell}
$\varphi_2 P_\hbar$ is elliptic on $\operatorname{WF}_\hbar( \varphi_1\Pi_H \varphi)$.
\ele
\begin{proof}
By (\ref{eq:WFprod}) and (\ref{eq:support}),
$\operatorname{WF}_\hbar (\varphi_1\Pi_H \varphi) \subset
\overline{T}^*(\operatorname{supp}\varphi_1) \cap \operatorname{WF}_\hbar \Pi_H$. Since $\varphi_2 = 1$ on $\operatorname{supp}\varphi_1$, the result is a direct consequence of Lemma \ref{lem:Pelliptic}.
\end{proof}

\

By the facts about compactly-supported operators recalled in \S\ref{sec:recap}, $\varphi_1\Pi_H \varphi$ and $\varphi_2 P_\hbar$ are compactly supported. Therefore, by Lemma \ref{lem:rev_ell}, we can apply Part (ii) of  Theorem \ref{thm:elliptic} with $A = \varphi_1\Pi_H \varphi$, $B_1 = \varphi_2$, $P = P_\hbar$, $m_1 = 0$, $m_2 = 0$, $\ell = 2$. This result implies that there exists $\widetilde \chi \in C^\infty_{\rm comp}$, and, for any $N'\geq 1$, there exists $C_{N'}>0$ such that
\beq \label{eq:rev_3}
\Vert \varphi_1\Pi_H \varphi u \Vert_{H^2_\hbar} \lesssim \Vert \varphi_2 P_\hbar u \Vert_{L^2} + C_{N'} \hbar^{N'} \Vert \widetilde \chi u \Vert_{L^2} = \hbar ^ 2\Vert \varphi_2 f \Vert_{L^2} + C_{N'} \hbar^{N'} \Vert \widetilde \chi u \Vert_{L^2}.
\eeq
Collecting (\ref{eq:rev_1}),  (\ref{eq:rev_2}),  (\ref{eq:rev_3}), using \eqref{eq:res}, and choosing $N=N'= M+1$, we obtain \eqref{eq:HF2}.

\section{Proof of Theorem \ref{thm:main2}}\label{sec:proof_FEM}

The two ingredients for the proof of Theorem \ref{thm:main2} are 
\bit
\item Lemma \ref{lem:Schatz}, which is the standard duality argument giving a condition for quasi-optimality to hold in terms of how well the solution of the adjoint problem 
is approximated by the finite-element space (measured by the quantity  $\eta(V_N)$ defined by \eqref{eq:etadef}), and  
\item Lemma \ref{lem:eta} that bounds  $\eta(V_N)$ using the decomposition from Theorem \ref{thm:main1}.
\eit 
Regarding Lemma \ref{lem:Schatz}: we recall that this argument came out of ideas introduced in \cite{Sc:74}, was then formalised in \cite{Sa:06}, and has been used extensively in the analysis of the Helmholtz FEM; see, e.g., \cite{AzKeSt:88, IhBa:95a, Me:95, Sa:06, MeSa:10, MeSa:11, ZhWu:13, Wu:14, DuWu:15, ChNi:18, LiWu:19, ChNi:20, GaChNiTo:18, GrSa:20, GaSpWu:20}.

Before stating Lemma \ref{lem:Schatz} we need to introduce some notation.

\begin{definition}[The adjoint sesquilinear form $a^*(\cdot,\cdot)$]\label{def:adjoint}
The adjoint sesquilinear form, $a^*(u,v)$, to the sesquilinear form $a(\cdot,\cdot)$ defined in \eqref{eq:sesqui} is given by
\beqs
a^*(u,v) := \overline{a(v,u)}= \int_{B_R} 
\Big((\MA \gu)\cdot\gvb
 - k^2 n u\vb\Big) - \big\langle \gamma u, \DtN(\gamma v)\big\rangle_{\partial B_R}.
\eeqs
\end{definition}

A key role is played by the solution operator of the adjoint variational problem  with data in $L^2(B_R)$; we therefore introduce the following notation.

\begin{definition}[Adjoint solution operator $\cS^*$]
Given $f\in L^2(B_R)$, let $\cS^*f$ be defined as the solution of the variational problem
\beq\label{eq:S*vp}
\tfind \cS^*f \in H^1(B_R) \quad\tst\quad a^*(\cS^*f,v) = \int_{B_R} f\, \overline{v} \quad\tfa v\in H^1(B_R).
\eeq
\end{definition}

Green's second identity applied to solutions of the Helmholtz equation satisfying the Sommerfeld radiation condition \eqref{eq:src} implies that
$\big\langle \DtN \psi, \overline{\phi}\big\rangle_{\GR} =\big\langle \DtN \phi, \overline{\psi}\big\rangle_{\GR} $
(see, e.g., \cite[Lemma 6.13]{Sp:15}); thus $a(\overline{v},u) = a(\overline{u},v)$ and so the definition \eqref{eq:S*vp} implies that
\beq\label{eq:S*fkey}
a(\overline{\cS^*f}, v)= (\overline{f},v)_{L^2(\OR)}\quad\tfa v\in H^1(B_R).
\eeq

\begin{definition}[$\eta(V_N)$]
Given a sequence of finite-dimensional spaces $(V_N)_{N=0}^\infty$ (as described in \S\ref{sec:FEMsetup}),
let
\beq\label{eq:etadef}
\eta(V_N):= \sup_{0\neq f\in L^2(B_R)}\min_{v_N\in V_N} \frac{\N{S^*f-v_N}_{H^1_k(B_R)}}{\big\|
f\big\|_{L^2(B_R)}}.
\eeq
\end{definition} 

\ble[Conditions for quasi-optimality]\label{lem:Schatz}
If 
\beqs
k\,\eta(V_N) \leq \frac{1}{\Ccont} \sqrt{ \frac{A_{\min}}{2\big(n_{\max} + A_{\min}\big)}},
\eeqs
then the Galerkin equations \eqref{eq:FEM} have a unique solution which satisfies
\beqs
\N{u-u_\hFEM}_{H^1_k(B_R)} \leq \frac{2\Ccont}{A_{\min}}\left(\min_{v_N\in V_N} \N{u-v_N}_{H^1_k(B_R)}\right).
\eeqs
\ele

\bpf
Using the inequality \eqref{eq:CDTN2}, we see that
$a(\cdot,\cdot)$ satisfies the G\aa rding inequality
\beq\label{eq:Garding}
\Re\big(a(v,v)\big) \geq A_{\min} \N{v}^2_{H^1_k(B_R)} - 2 k^2 \big(n_{\max} + A_{\min}\big)\N{v}^2_{L^2(B_R)}
\eeq
and the result follows from, e.g., the account \cite[Theorem 6.32]{Sp:15} of the standard duality argument with (in the notation of \cite{Sp:15}) $\alpha = A_{\min}$ and $C_{\mathcal{V}}= 2 k^2 \big(n_{\max} + A_{\min}\big)$.
\epf

\ble[Bound on $\eta(V_N)$ using the decomposition from Theorem \ref{thm:main1}]\label{lem:eta}
Let $\MA$ and $n$ satisfy Assumption \ref{ass:1} and let $R>0$ be such that $\supp(\MI-\MA)\cup\supp(1-n)\Subset B_R$. 
Let $(V_N)_{N=0}^\infty$ be the piecewise-polynomial approximation spaces described in \cite[\S5]{MeSa:10}. 
There exists $C_6,C_7, \sigma>0$, all independent of $k,h,$ and $p$, such that 
\beq\label{eq:etabound}
k\,\eta(V_N) \leq 
 C_6 C_3\frac{hk}{p}\left(1 + \frac{kh}{p}\right) 
 +
C_7 \Csol(k) \left[
\left(\frac{h}{h+\sigma}\right)^p \left(1 + \frac{hk}{h+\sigma}\right) + k \left(\frac{kh}{\sigma p}\right)^p 
\left(\frac{1}{p} + \frac{kh}{\sigma p}\right) 
 \right].
\eeq
The constants $C_6$ and $\sigma$ only depend on the constants in \cite[Assumption 5.2]{MeSa:10} defining the element maps from the reference element;
$C_7$ depends on these constants, and additionally on $C_5$.
\ele

\bpf
This proof is very similar to the proof of \cite[Theorem 5.5]{MeSa:10}. Indeed, \cite[Theorem 5.5]{MeSa:10} proves a bound very similar to \eqref{eq:etabound} starting from bounds almost identical to the bounds \eqref{eq:uhigh} and \eqref{eq:ulow} (recalling Remark \ref{rem:nabla} about notation). The only difference is that the bound \eqref{eq:ulow} contains $\Csol$, which depends on $k$ (whereas in \cite{MeSa:10} $\Csol \sim 1$), and so we now need to keep track of how $\Csol$ enters the proof of \cite[Theorem 5.5]{MeSa:10}.

From the definition \eqref{eq:etadef}, it is sufficient to show that, given $f\in L^2(B_R)$, there exists $w_N\in V_N$ such that
\beq\label{eq:Sat1}
\N{\mathcal{S}^* f - w_N}_{H^1_k(B_R)}\leq C \N{f}_{L^2(B_R)},
\eeq
where $C$ is the right-hand side of \eqref{eq:etabound} divided by $k$.
Let $v:= \mathcal{S}^* f$; by \eqref{eq:S*fkey} and Part (i) of Lemma \ref{lem:BVPs}, $\overline{v}$ satisfies the assumptions of Theorem \ref{thm:main1} with $f$ replaced by $\overline{f}$, and so the bounds \eqref{eq:uhigh} and \eqref{eq:ulow} hold with $u$ replaced by $v$.

By \cite[First equation on Page 1896]{MeSa:10} (which uses \cite[Theorem B.4]{MeSa:10}), 
the bound \eqref{eq:MSbae} holds, and thus 
there exists 
$w_N^{(1)} \in V_N$ such that 
\beqs
\N{\vhigh - w_N^{(1)}}_{H^1_k(B_R)} \leq C_6 \frac{h}{p}\left(1 + \frac{kh}{p}\right)  |v|_{H^2(B_R)}
\eeqs
and so 
\beq\label{eq:uhighapprox}
\N{\vhigh - w_N^{(1)}}_{H^1_k(B_R)} \leq C_6 \frac{h}{p}\left(1 + \frac{kh}{p}\right) C_3\N{f}_{L^2(B_R)}
\eeq
by \eqref{eq:uhigh}.

For the approximation of $\vlow$, the only change to the argument in \cite{MeSa:10} is that a multiplicative factor of $(\Csol)^2$ must be included on the right-hand side of \cite[Equation 5.8]{MeSa:10}. Then \cite[Equations 5.8 and 5.9]{MeSa:10} implies that 
there exists $C_7$ and $w_N^{(2)} \in V_N$ such that 
\beq\label{eq:ulowapprox}
k \N{\vlow - w_N^{(2)}}_{H^1_k(B_R)} \leq C_7 \Csol(k) \left[
\left(\frac{h}{h+\sigma}\right)^p \left(1 + \frac{hk}{h+\sigma}\right) + k \left(\frac{kh}{\sigma p}\right)^p 
\left(\frac{1}{p} + \frac{kh}{\sigma p}\right) 
 \right]\N{f}_{L^2(B_R)}
\eeq
(observe that this equation is identical to \cite[Last equation on Page 1896]{MeSa:10} except for the factor $\Csol$ on the right-hand side).

Let $w_N:= w_N^{(1)}+ w_N^{(2)}$. By the triangle inequality, the decomposition $v=\vhigh+ \vlow$ on $B_R$, and the inequalities \eqref{eq:uhighapprox} and \eqref{eq:ulowapprox}, the inequality \eqref{eq:Sat1} holds with $C$ the right-hand side of \eqref{eq:etabound} and the proof is complete.
\epf

\begin{corollary}[Conditions under which $k\,\eta(V_N)$ is arbitrarily small]\label{cor:etasmall}
Let the assumptions of Lemma \ref{lem:eta} hold. Given $\eps>0$ and $k_0>0$, there exists $\mathcal{C}_1, \mathcal{C}_2>0$, depending only on $\eps, C_3, C_6, C_7, \sigma,$ and $k_0$,
such that if 
\beqs
\frac{hk}{p}\leq \mathcal{C}_1 \quad\tand\quad p \geq \mathcal{C}_2 \Big(1 + \log k + \log \big(\Csol(k)\big)\Big),
\eeqs
then
\beqs
k\, \eta(V_N)\leq \eps \quad\tfa k\geq k_0.
\eeqs
\end{corollary}

\bpf
This proof is essentially identical to the proofs of \cite[Corollary 5.6]{MeSa:10} and \cite[Theorem 5.8]{MeSa:11}. 
First choose $\mathcal{C}_1$ sufficiently small such that $\mathcal{C}_1< \sigma$ and 
\beqs
 C_6 \,C_3\,
 \mathcal{C}_1
 \left(1 + 
 \mathcal{C}_1
 \right) \leq \frac{\eps}{2}
\eeqs
From the bound on $k\eta(V_N)$ \eqref{eq:etabound}, it is then sufficient to show that 
\beq\label{eq:small}
C_7\, \Csol(k) \left[
\left(\frac{h}{h+\sigma}\right)^p \left(1 + \frac{hk}{h+\sigma}\right) + k \left(\frac{kh}{\sigma p}\right)^p 
\left(\frac{1}{p} + \frac{kh}{\sigma p}\right) 
 \right]
\eeq
can be made $\leq \eps/2$. Let 
\beqs
\theta_1: = \frac{h}{h+\sigma} \quad\tand\quad \theta_2:= \frac{\mathcal{C}_1}{\sigma},
\eeqs
so that  \eqref{eq:small} is bounded by
\beqs
C_7 \,\Csol(k) \left[ (\theta_1)^p \left(1+ \frac{\mathcal{C}_1 p }{\sigma}\right) + k (\theta_2)^p \left(\frac{1}{p}+ \frac{\mathcal{C}_1}{\sigma}\right)\right];
\eeqs
the result then follows since $\theta_1, \theta_2 <1$.
\epf

\

\bpf[Proof of Theorem \ref{thm:main2}]
This follows by combining Lemma \ref{lem:Schatz} and Corollary \ref{cor:etasmall}.
\epf

\section*{Acknowledgements}

The authors thank 
Martin Averseng (ETH Z\"urich) and an anonymous referee for highlighting simplifications of the arguments in a earlier version of the paper. We also thank Th\'eophile Chaumont-Frelet (INRIA, Nice) for useful
discussions about the results of \cite{MeSa:10}, \cite{MeSa:11}. 
DL and EAS acknowledge support from EPSRC grant EP/1025995/1. JW was partly supported by Simons Foundation grant 631302.

\footnotesize{
\bibliographystyle{siam}
\bibliography{biblio_combined_sncwadditions}

\begin{thebibliography}{10}

\bibitem{AzKeSt:88}
{\sc A.~K. Aziz, R.~B. Kellogg, and A.~B. Stephens}, {\em A two point boundary
  value problem with a rapidly oscillating solution}, Numer. Math., 53 (1988),
  pp.~107--121.

\bibitem{BaSa:00}
{\sc I.~M. Babu\v{s}ka and S.~A. Sauter}, {\em Is the pollution effect of the
  {FEM} avoidable for the {H}elmholtz equation considering high wave numbers?},
  SIAM Review,  (2000), pp.~451--484.

\bibitem{BaChGo:17}
{\sc H.~Barucq, T.~Chaumont-Frelet, and C.~Gout}, {\em {Stability analysis of
  heterogeneous Helmholtz problems and finite element solution based on
  propagation media approximation}}, Math. Comp., 86 (2017), pp.~2129--2157.

\bibitem{Me:20}
{\sc M.~Bernkopf, T.~Chaumont-Frelet, and J.~M. Melenk}, {\em {Stability and
  convergence of Galerkin discretizations of the Helmholtz equation in
  piecewise smooth media}},
  \url{https://numericalwaves.sciencesconf.org/data/program/abstract_melenk.pdf},
   (2020).

\bibitem{Bu:98}
{\sc N.~Burq}, {\em D\'ecroissance des ondes absence de de l'\'energie locale
  de l'\'equation pour le probl\`{e}me ext\'erieur et absence de resonance au
  voisinage du r\'eel}, Acta Math., 180 (1998), pp.~1--29.

\bibitem{Bu:02}
\leavevmode\vrule height 2pt depth -1.6pt width 23pt, {\em Semi-classical
  estimates for the resolvent in nontrapping geometries}, International
  Mathematics Research Notices, 2002 (2002), pp.~221--241.

\bibitem{ChMo:08}
{\sc S.~N. Chandler-Wilde and P.~Monk}, {\em {Wave-number-explicit bounds in
  time-harmonic scattering}}, SIAM J. Math. Anal., 39 (2008), pp.~1428--1455.

\bibitem{Ch:16}
{\sc T.~Chaumont-Frelet}, {\em {On high order methods for the heterogeneous
  Helmholtz equation}}, Computers \& Mathematics with Applications, 72 (2016),
  pp.~2203--2225.

\bibitem{ChNi:18}
{\sc T.~Chaumont-Frelet and S.~Nicaise}, {\em {High-frequency behaviour of
  corner singularities in Helmholtz problems}}, ESAIM: Math. Model. Numer.
  Anal., 52 (2018), pp.~1803--1845.

\bibitem{ChNi:20}
\leavevmode\vrule height 2pt depth -1.6pt width 23pt, {\em {Wavenumber explicit
  convergence analysis for finite element discretizations of general wave
  propagation problem}}, IMA J. Numer. Anal., 40 (2020), pp.~1503--1543.

\bibitem{CoDaNi:10}
{\sc M.~Costabel, M.~Dauge, and S.~Nicaise}, {\em {Corner Singularities and
  Analytic Regularity for Linear Elliptic Systems. Part I: Smooth domains.}},
  (2010).
\newblock
  \url{https://hal.archives-ouvertes.fr/file/index/docid/453934/filename/CoDaNi_Analytic_Part_I.pdf}.

\bibitem{DuWu:15}
{\sc Y.~Du and H.~Wu}, {\em {Preasymptotic error analysis of higher order FEM
  and CIP-FEM for Helmholtz equation with high wave number}}, SIAM J. Numer.
  Anal., 53 (2015), pp.~782--804.

\bibitem{DuHo:72}
{\sc J.~J. Duistermaat and L.~H{\"o}rmander}, {\em Fourier integral operators.
  ii}, Acta mathematica, 128 (1972), pp.~183--269.

\bibitem{DyZw:19}
{\sc S.~Dyatlov and M.~Zworski}, {\em Mathematical theory of scattering
  resonances}, vol.~200 of Graduate Studies in Mathematics, American
  Mathematical Society, 2019.

\bibitem{EsMe:12}
{\sc S.~Esterhazy and J.~M. Melenk}, {\em On stability of discretizations of
  the {H}elmholtz equation}, in Numerical Analysis of Multiscale Problems,
  I.~G. Graham, T.~Y. Hou, O.~Lakkis, and R.~Scheichl, eds., Springer, 2012,
  pp.~285--324.

\bibitem{GaSpWu:20}
{\sc J.~Galkowski, E.~A. Spence, and J.~Wunsch}, {\em {Optimal constants in
  nontrapping resolvent estimates}}, Pure and Applied Analysis, 2 (2020),
  pp.~157--202.

\bibitem{GaChNiTo:18}
{\sc D.~Gallistl, T.~Chaumont-Frelet, S.~Nicaise, and J.~Tomezyk}, {\em
  Wavenumber explicit convergence analysis for finite element discretizations
  of time-harmonic wave propagation problems with perfectly matched layers},
  hal preprint 01887267,  (2018).

\bibitem{GaMo:19}
{\sc M.~Ganesh and C.~Morgenstern}, {\em {A coercive heterogeneous media
  Helmholtz model: formulation, wavenumber-explicit analysis, and
  preconditioned high-order FEM}}, Numerical Algorithms,  (2019), pp.~1--47.

\bibitem{GoGrSp:20}
{\sc S.~Gong, I.~G. Graham, and E.~A. Spence}, {\em {Domain decomposition
  preconditioners for high-order discretisations of the heterogeneous Helmholtz
  equation}}, IMA J. Num. Anal., 41 (2021), pp.~2139--2185.

\bibitem{GrPeSp:19}
{\sc I.~G. Graham, O.~R. Pembery, and E.~A. Spence}, {\em {The Helmholtz
  equation in heterogeneous media: a priori bounds, well-posedness, and
  resonances}}, Journal of Differential Equations, 266 (2019), pp.~2869--2923.

\bibitem{GrSa:20}
{\sc I.~G. Graham and S.~A. Sauter}, {\em {Stability and finite element error
  analysis for the Helmholtz equation with variable coefficients}}, Math.
  Comp., 89 (2020), pp.~105--138.

\bibitem{HeRo:83}
{\sc B.~Helffer and D.~Robert}, {\em Calcul fonctionnel par la transformation
  de {M}ellin et op\'{e}rateurs admissibles}, J. Funct. Anal., 53 (1983),
  pp.~246--268.

\bibitem{HeSj:89}
{\sc B.~Helffer and J.~Sj\"{o}strand}, {\em \'{E}quation de {S}chr\"{o}dinger
  avec champ magn\'{e}tique et \'{e}quation de {H}arper}, in Schr\"{o}dinger
  operators ({S}\o nderborg, 1988), vol.~345 of Lecture Notes in Phys.,
  Springer, Berlin, 1989, pp.~118--197.

\bibitem{Ho:83}
{\sc L.~H\"{o}rmander}, {\em The Analysis of Linear Differential Operators. I,
  Distribution Theory and Fourier Analysis}, Springer-Verlag, Berlin, 1983.

\bibitem{HsWe:08}
{\sc G.~C. Hsiao and W.~L. Wendland}, {\em Boundary integral equations},
  vol.~164 of Applied Mathematical Sciences, Springer, 2008.

\bibitem{IhBa:95a}
{\sc F.~Ihlenburg and I.~Babu{\v{s}}ka}, {\em {Finite element solution of the
  Helmholtz equation with high wave number Part I: The h-version of the FEM}},
  Comput. Math. Appl., 30 (1995), pp.~9--37.

\bibitem{LaSpWu:21}
{\sc D.~Lafontaine, E.~A. Spence, and J.~Wunsch}, {\em {Decompositions of
  high-frequency Helmholtz solutions via functional calculus, and application
  to the finite element method}}, arXiv preprint arXiv:2102.13081,  (2021).

\bibitem{LaSpWu:20}
\leavevmode\vrule height 2pt depth -1.6pt width 23pt, {\em For most
  frequencies, strong trapping has a weak effect in frequency-domain
  scattering}, Communications on Pure and Applied Mathematics, 74 (2021),
  pp.~2025--2063.

\bibitem{LaSpWu:19}
\leavevmode\vrule height 2pt depth -1.6pt width 23pt, {\em {A sharp
  relative-error bound for the Helmholtz $h$-FEM at high frequency}},
  Numerische Mathematik, 150 (2022), pp.~137--178.

\bibitem{LaxPhi}
{\sc P.~D. Lax and R.~S. Phillips}, {\em {Scattering Theory}}, Academic Press,
  revised~ed., 1989.

\bibitem{LiWu:19}
{\sc Y.~Li and H.~Wu}, {\em {FEM and CIP-FEM for Helmholtz Equation with High
  Wave Number and Perfectly Matched Layer Truncation}}, SIAM J. Numer. Anal.,
  57 (2019), pp.~96--126.

\bibitem{Ma:02a}
{\sc A.~Martinez}, {\em An introduction to semiclassical and microlocal
  analysis}, vol.~994, Springer, 2002.

\bibitem{Me:95}
{\sc J.~M. Melenk}, {\em {On generalized finite element methods}}, PhD thesis,
  The University of Maryland, 1995.

\bibitem{MePaSa:13}
{\sc J.~M. Melenk, A.~Parsania, and S.~Sauter}, {\em {General DG-methods for
  highly indefinite Helmholtz problems}}, Journal of Scientific Computing, 57
  (2013), pp.~536--581.

\bibitem{MeSa:10}
{\sc J.~M. Melenk and S.~Sauter}, {\em Convergence analysis for finite element
  discretizations of the {H}elmholtz equation with {D}irichlet-to-{N}eumann
  boundary conditions}, Math. Comp, 79 (2010), pp.~1871--1914.

\bibitem{MeSa:11}
\leavevmode\vrule height 2pt depth -1.6pt width 23pt, {\em Wavenumber explicit
  convergence analysis for {G}alerkin discretizations of the {H}elmholtz
  equation}, SIAM J. Numer. Anal., 49 (2011), pp.~1210--1243.

\bibitem{Ne:01}
{\sc J.~C. N{\'e}d{\'e}lec}, {\em {Acoustic and electromagnetic equations:
  integral representations for harmonic problems}}, Springer Verlag, 2001.

\bibitem{Pe:20}
{\sc O.~R. Pembery}, {\em {The Helmholtz Equation in Heterogeneous and Random
  Media: Analysis and Numerics}}, PhD thesis, University of Bath, 2020.

\bibitem{Ra:71}
{\sc J.~V. Ralston}, {\em Trapped rays in spherically symmetric media and poles
  of the scattering matrix}, Communications on Pure and Applied Mathematics, 24
  (1971), pp.~571--582.

\bibitem{Rob87}
{\sc D.~Robert}, {\em Autour de l'approximation semi-classique}, vol.~68 of
  Progress in Mathematics, Birkh\"{a}user Boston, Inc., Boston, MA, 1987.

\bibitem{SaVa:02}
{\sc J.~Saranen and G.~Vainikko}, {\em Periodic integral and pseudodifferential
  equations with numerical approximation}, Springer, 2002.

\bibitem{Sa:06}
{\sc S.~A. Sauter}, {\em {A refined finite element convergence theory for
  highly indefinite Helmholtz problems}}, Computing, 78 (2006), pp.~101--115.

\bibitem{Sc:74}
{\sc A.~H. Schatz}, {\em {An observation concerning Ritz-Galerkin methods with
  indefinite bilinear forms}}, Math. Comp., 28 (1974), pp.~959--962.

\bibitem{Sj:97}
{\sc J.~Sj\"{o}strand}, {\em A trace formula and review of some estimates for
  resonances}, in Microlocal analysis and spectral theory ({L}ucca, 1996),
  vol.~490 of NATO Adv. Sci. Inst. Ser. C Math. Phys. Sci., Kluwer Acad. Publ.,
  Dordrecht, 1997, pp.~377--437.

\bibitem{SjZw:91}
{\sc J.~Sj\"{o}strand and M.~Zworski}, {\em Complex scaling and the
  distribution of scattering poles}, J. Amer. Math. Soc., 4 (1991),
  pp.~729--769.

\bibitem{Sp:15}
{\sc E.~A. Spence}, {\em {Overview of Variational Formulations for Linear
  Elliptic PDEs}}, in Unified transform method for boundary value problems:
  applications and advances, A.~S. Fokas and B.~Pelloni, eds., SIAM, 2015,
  pp.~93--159.

\bibitem{Ta:96}
{\sc M.~E. Taylor}, {\em Partial differential equations II, Qualitative studies
  of linear equations, volume 116 of Applied Mathematical Sciences},
  Springer-Verlag, New York, 1996.

\bibitem{Va:75}
{\sc B.~R. Vainberg}, {\em {On the short wave asymptotic behaviour of solutions
  of stationary problems and the asymptotic behaviour as $t\rightarrow \infty$
  of solutions of non-stationary problems}}, Russian Mathematical Surveys, 30
  (1975), pp.~1--58.

\bibitem{Va:89}
\leavevmode\vrule height 2pt depth -1.6pt width 23pt, {\em Asymptotic methods
  in equations of mathematical physics}, Gordon \& Breach Science Publishers,
  New York, 1989.
\newblock Translated from the Russian by E. Primrose.

\bibitem{Wu:14}
{\sc H.~Wu}, {\em {Pre-asymptotic error analysis of CIP-FEM and FEM for the
  Helmholtz equation with high wave number. Part I: linear version}}, IMA J.
  Numer. Anal., 34 (2014), pp.~1266--1288.

\bibitem{ZhWu:13}
{\sc L.~Zhu and H.~Wu}, {\em {Preasymptotic error analysis of CIP-FEM and FEM
  for Helmholtz equation with high wave number. Part II: $hp$ version}}, SIAM
  J. Numer. Anal., 51 (2013), pp.~1828--1852.

\bibitem{Zw:12}
{\sc M.~Zworski}, {\em Semiclassical analysis}, vol.~138 of Graduate Studies in
  Mathematics, American Mathematical Society, Providence, RI, 2012.

\end{thebibliography}
}

\end{document}